\documentclass[10pt]{article}
\usepackage{epsfig}
\usepackage{amsmath,amssymb, amsfonts}
\usepackage{xspace}
\def\be{\begin{equation}}
\def\ee{\end{equation}}
\def\bea{\begin{eqnarray}}
\def\eea{\end{eqnarray}}
\def\bes{\begin{eqnarray*}}
\def\ees{\end{eqnarray*}}

\def\nn{\nonumber}
\def\<{\langle}
\def\>{\rangle}
\def\lb{\label}
\def\bs{\setminus}

\def\R{{\mathbb{R}}}
\def\C{{\mathbb{C}}}
\def\Z{{\mathbb{Z}}}
\def\N{{\mathbb{N}}}
\def\U{{\mathbb{U}}}
\def\Q{{\mathbb{Q}}}
\def\T{{\bf T}}

\def\RP{{\mathbb{R}P}}

\def\aa{{\alpha}}
\def\bb{{\beta}}
\def\ga{{\gamma}}

\def\th{{\theta}}
\def\om{{\omega}}
\def\Om{{\Omega}}
\def\ep{{\epsilon}}
\def\lm{{\lambda}}
\def\Lm{{\Lambda}}

\def\sg{{\sigma}}

\def\vf{{\varphi}}

\def\vs{{\varsigma}}

\def\P{{\cal P}}

\def\id{{\it id}}
\def\im{{\rm im}}

\def\rank{{\rm rank}}
\def\Sp{{\rm Sp}}

\def\ol#1{\overline{#1}}  %overline in math mode
\def\td#1{\tilde{#1}}
\def\hb{\vrule height0.18cm width0.14cm $\,$}

\def\ol#1{\overline{#1}}  %overline in math mode
\def\td#1{\tilde{#1}}

\title{The second closed geodesic on Finsler manifolds with $\Z_2$ fundamental group}
\author{Huagui Duan $^{1}$,\thanks{Partially supported by NSFC (12631007, 12271268 and 12361141812) and Natural Science
Foundation of Tianjin (No. 25JCZDJC01030). E-mail: duanhg@nankai.edu.cn.}
\quad Hui Liu $^{2}$,\thanks{Partially supported by NSFC (12371195 and 12022111) and the Fundamental
Research Funds for the Central Universities (2042023kf0207). E-mail: huiliu00031514@whu.edu.cn.}
\quad Wei Wang $^{3}$\thanks{Partially
supported by NSFC (12631007, 12025101 and 11431001). E-mail: wangwei@math.pku.edu.cn.}\\
$^{1}$ School of Mathematics, Jilin University, Changchun, 130012, Jilin\\
$^{2}$ School of Mathematics and Statistics, Wuhan University, Wuhan 430072, Hubei\\
$^{3}$ KLPAM and School of Mathematical Sciences, Peking University, Beijing 100871 \\
The People's Republic of China
}

\begin{document}
\newtheorem{definition}{Definition}[section]
\newtheorem{theorem}{Theorem}[section]
\newtheorem{lemma}{Lemma}[section]
\newtheorem{corollary}{Corollary}[section]
\newtheorem{example}{Example}[section]
\newtheorem{property}{Property}[section]
\newtheorem{proposition}{Proposition}[section]
\newtheorem{remark}{Remark}[section]

\newcommand{\qed}{\nolinebreak\hfill\rule{2mm}{2mm}}
\newcommand{\Proof}{\par\medbreak\it Proof: \rm}
\newcommand{\rem}{\par\medbreak\it Remark: \rm}
\newcommand{\defi}{\par\medbreak\it Definition : \rm}
\renewcommand{\thefootnote}{\arabic{footnote}}

\maketitle

\begin{abstract}
{\it  Let $(M, F)$ be a compact bumpy Finsler manifold with $\pi_1(M)=\Z_2$,
we prove that there exist at least two distinct non-contractible closed geodesics on $(M, F)$.}
\end{abstract}

{\bf Key words}: Non-contractible closed geodesics;  Non-simply connected manifolds; Morse theory; Index iteration theory;
Systems of irrational numbers; Kronecker's approximation theorem

{\bf AMS Subject Classification}: 53C22, 58E05, 58E10.

\renewcommand{\theequation}{\thesection.\arabic{equation}}
\makeatletter\@addtoreset{equation}{section}\makeatother

\setcounter{equation}{0}%\setcounter{figure}{0}
\section{Introduction and sketch of proofs}%S 1.

\subsection{Introduction }  There is a famous  conjecture in Riemannian geometry which claims the existence of
infinitely many  distinct closed geodesics on every compact Riemannian manifold $(M, g)$ with $\dim M\geq2$.
By Gromoll and Meyer's results  in \cite{GrM} of 1969,  if the Betti number sequence $\{b_p(\Lm M)\}_{p\in\N}$ of the free
loop space $\Lm M$ of $M$ is unbounded, then the conjecture is true. Then in \cite{VPS} of 1976, for compact simply
connected manifold $M$, Vigu\'e-Poirrier and Sullivan further proved this Betti number
sequence is bounded if and only if $M$ satisfies
\be  H^*(M;\Q)\cong T_{d,n+1}(x)=\Q[x]/(x^{n+1}=0)   \lb{1.1}\ee
with a generator $x$ of degree $d\ge 2$ and height $n+1\ge 2$, where $\dim M=dn$.
Typical examples of these manifolds are CROSSs (compact rank one symmetric spaces),
e.g., $S^d$, complex projective spaces  $\C P^n$ (with $d=2$), quaternionic
projective spaces $H P^n$ (with $d=4$), and the Cayley plane $CaP^2$ (with $d=8$).
Among these manifolds, only for Riemannian $S^2$ this conjecture was proved by
the works \cite{Fra} of Franks in 1992 and \cite{Ban} of Bangert in 1993.

On the other hand, this conjecture is false for  Finsler
manifolds. In 1973,  Katok in \cite{Kat} found some irreversible Finsler metrics on CROSSs which possess
precisely $\frac{dn(n+1)}{2}$ ($d+1$) distinct closed geodesics for even (odd) $d$  and all of which are non-degenerate and elliptic.

On a Finsler
manifold $(M,F)$, the $m$-th iterate $c^m$ of a closed curve $c:S^1=\R/\Z\to M$ is defined
by $c^m(t)=c(mt)$. A closed geodesic $c$ is {\it prime} if it is not a multiple covering
(i.e., iterate) of any other closed geodesics. The inverse curve $c^{-1}$ of $c$ is defined by
$c^{-1}(t)=c(1-t)$ for $t\in \R$. We call two prime closed geodesics $c$ and $d$ {\it distinct}, if there is no $\th\in [0,1)$
such that $c(t)=d(t+\th)$ for all $t\in\R$. We shall omit the word {\it distinct} when we talk
about more than one prime closed geodesic. On a reversible Finsler (or Riemannian) manifold,
two prime closed geodesics $c$ and $d$ are called { \it  distinct} if $c(S^1)\neq d(S^1)$
in $M$. Equivalently, there exist no \(\varepsilon\in\{1,-1\}\) and \(\theta\in[0,1)\) such that \(d(t)=c(\varepsilon t+\theta)\) for all \(t\in\mathbb R\).

For a closed geodesic $c$ on $(M,F)$ with $n=\dim M$, denote by $P_c$ the linearized
Poincar\'{e} map of $c$. Then $P_c\in \Sp(2n-2)$ is symplectic. A closed geodesic $c$ is called
{\it non-degenerate} if $1$ is not an eigenvalue of $P_c$. A Finsler manifold $(M,F)$ is called
{\it bumpy} if all the closed geodesics on $M$ are non-degenerate.

By the classical theorem of Lyusternik-Fet \cite{LyF} in 1951, there exists at least one closed geodesic on every compact
Riemannian manifold, their proof is variational and works also for the Finsler manifolds.
In 1965, A. Fet proved in \cite{Fet}
that for compact simply connected bumpy Riemannian (or reversible Finsler) manifolds, the existence of at least two distinct closed geodesics, while his proof relies on the symmetric property
of the metric, hence can not be applied to  irreversible Finsler metrics.
In \cite{DLW1}, H. Duan, Y. Long and W. Wang proved   the existence of at least two distinct closed geodesics on any compact simply connected bumpy  irreversible Finsler manifolds.

The following are the main results in this paper.

{\bf Theorem 1.1.} {\it Suppose $(M, F)$ be a compact bumpy irreversible Finsler manifold satisfying $\pi_1(M)=\Z_2$.
Then there are at least two distinct non-contractible closed geodesics on $(M, F)$.}

{\bf Theorem 1.2.} {\it Suppose $(M, F)$ be a compact bumpy reversible Finsler manifold (or Riemannian manifold) satisfying $\pi_1(M)=\Z_2$.
Then there are at least two distinct non-contractible closed geodesics on $(M, F)$.}

We are aware of a number of results concerning closed geodesics on Finsler manifolds.  In \cite{BaL}, V. Bangert and Y. Long proved that on any
Finsler 2-sphere $(S^2, F)$, there exist at least two  closed
geodesics. In \cite{Rad5}, H.-B. Rademacher studied the existence and
stability of closed geodesics on positively curved Finsler
manifolds. In \cite{Wan1}-\cite{Wan3}, W. Wang studied the existence and
stability of closed geodesics on positively curved Finsler spheres.
In \cite{DuL1} of H. Duan and Y. Long and in  \cite{Rad6} of H.-B. Rademacher, they proved there exist
at least two  closed geodesics on any bumpy
Finsler $n$-sphere independently. In  \cite{Rad7}, H.-B. Rademacher proved there exist at least two  closed geodesics on any bumpy Finsler $\mathbb{C}P^2$. In \cite{DLW2}, H. Duan, Y. Long and W. Wang proved there exist at least $\frac{dn(n+1)}{2}$
non-hyperbolic closed geodesics on any compact simply-connected bumpy Finsler manifold $(M,F)$  with $H^*(M;\Q)\cong T_{d,n+1}(x)$ and flag curvature $K\ge 0$, where $d\geq2$ is even. In \cite{LiX}, H. Liu and Y. Xiao proved there exist at least two non-contractible  closed geodesics on
any bumpy Finsler ${ \RP}^n$.
In \cite{LLX}, H. Liu, Y. Long and Y. Xiao
proved there exist at least two non-contractible  closed geodesics on any bumpy Finsler $S^n/\Gamma$,
where $\Gamma$ is a finite group acts on $S^n$ freely  and isometrically.
In \cite{Liu}, H. Liu proved there exist $2[\frac{n+1}{2}]$  closed geodesics on any bumpy Finsler $S^n/\Gamma$ under a pinching condition.  In \cite{GG} and \cite{GGM} V. Ginzburg-B. Gurel and V. Ginzburg-B. Gurel-L. Macarini obtains some lower bounds for the number of closed geodesics
on certain Finsler manifolds under weak index conditions.

{\bf Remark 1.3.} Note that  Theorems 1.1 and 1.2 can not be deduced from
the results in \cite{DLW1}. In fact, one can find two closed geodesics $\tilde c$ and $\tilde d$
on the universal covering space $\tilde M$ of $M$, while one can not obtain that
whether $\tilde c$ or $\tilde d$ comes from the iteration of two distinct non-contractible
closed geodesics $c$ and $d$, i.e., $\tilde c=c^2$ and $\tilde d=d^2$.

\subsection{Sketch of proofs of Theorem 1.1}

Let $\tilde M$ be the universal covering space of $M$ and $g: \tilde M\rightarrow \tilde M$
be the deck transformation of $g\in\pi_1(M)$, where $\pi_1(M)=\{e, g\}$ with $e$ being
the identity and $g$ being the generator of $\Z_2$ satisfying $g^2=e$. Then the free loop space $\Lambda M$ possesses
a natural decomposition\bea \Lambda M=\Lambda_e M\bigsqcup\Lambda_g M,\lb{1.2}\eea
where $\Lambda_e M$ and $\Lambda_g M$ are the two connected components of $\Lambda M$ whose elements are homotopic
to $e$ and $g$ respectively.  Thus Theorem 1.1 says that we  can find at least two distinct closed geodesics in $\Lambda_g M$.

Let $X=\tilde M$ and $G(g)=\{(x,g(x))\mid x\in X\}$. Following \cite{GKV}, we consider the path space
\begin{equation}
P_gX=\tilde M_{G(g)}^I
=\{\gamma:\R\longrightarrow X\mid\gamma(t+1)=g\gamma(t)\}.
\label{1.4}
\end{equation}
The covering projection $p:X\longrightarrow M$ induces a double cover
\begin{equation}
\pi:P_gX\longrightarrow\Lambda_gM,\qquad \pi(\gamma)=p\circ\gamma,
\label{1.5}
\end{equation}
whose deck transformation is $\Gamma(\gamma)=g\gamma=\gamma(\cdot+1)$.
Thus $\Lambda_gM\cong P_gX/\langle\Gamma\rangle$. Since $\Gamma$ is
homotopic to the identity through parameter translations, the transfer gives
$H^*(\Lambda_gM;\Q)\cong H^*(P_gX;\Q)$.

Let $T=\R/2\Z$ and $\rho:T\longrightarrow S^1=\R/\Z$ be the natural homomorphism.
The parameter translation defines a $T$-action on $P_gX$, and $\pi$ is
$\rho$-equivariant. The induced map $B\rho$ is a rational cohomology
isomorphism. As in the proof of Proposition 3.2, we obtain
\bea
H^*_{S^1}(\Lambda_gM;\Q)\cong H^*_{T}(P_gX;\Q).
\eea
$g$ induces a homomorphism of homotopy groups
\bea g_\#:\pi_\ast (\tilde M)\rightarrow\pi_\ast (\tilde M).\lb{1.6}\eea
As in P.278 of \cite{GKV}, denote by $\pi_\ast(\tilde M)^{g_\#}$
the homotopy of $\tilde M$ fixed by $g_\#$.
Then $\{b_i(\tilde M_{G(g)}^I)\}_{i\in\N}$ is bounded if and only if
\bea \dim \pi^{even}_\ast(\tilde M)^{g_\#}\otimes \Q
\le \dim \pi^{odd}_\ast(\tilde M)^{g_\#}\otimes \Q\le 1.\lb{1.7}\eea

Choose a $g$-equivariant minimal Sullivan model
$\phi:A=(\Lambda V,d)\longrightarrow A_{PL}(X)$ with an involution $h$
satisfying $h^2=1$, $h(V)=V$ and $\phi h=g^*\phi$; see \cite{FHT}.
Write $V=V^+\oplus V^-$ for the eigenspace decomposition of $h$.
The ideal $(V^-)$ is $d$-stable, and hence
\bea
B=A/(V^-)=(\Lambda Y,\delta),\qquad Y=V^+,
\eea
is a simply connected minimal Sullivan algebra of finite type.
In the proof of Proposition 3.2, its Sullivan realization $Z$ and the
quotient map $A\longrightarrow B$ give equivariant comparisons, which yield
\bea
\begin{aligned}
H^*(\Lambda_gM;\Q)&\cong H^*(\Lambda Z;\Q),\\
H^*_{S^1}(\Lambda_gM;\Q)&\cong H^*_{S^1}(\Lambda Z;\Q).
\end{aligned}
\eea
Therefore, the $S^1$-equivariant cohomology is computed by
\bea
\mathcal C_B=\bigl(\Q[u]\otimes\Lambda(Y\oplus sY),D\bigr),
\qquad |u|=2,\quad |sv|=|v|-1,
\eea
where $Du=0$, $Dv=\delta v+u\,sv$ and $D(sv)=-s(\delta v)$;
see \cite{VPS,VPB}. Setting $u=0$ gives the ordinary free-loop model. Then the classification in Theorem 4.5 of \cite{GKV} takes the following
form.

{\bf Classification Theorem.} {\it Denote by $\chi_e=\dim Y^{even}$ and
$\chi_o=\dim Y^{odd}$. Exactly one of the following holds.

(i) $\{b_q(P_gX)\}_{q\in\N}$ is unbounded.

(ii) $\chi_o=\chi_e=0$. In this case $B=\Q$ and $H^*(P_gX;\Q)=\Q$.

(iii) $\chi_o=1$ and $\chi_e=0$. In this case
$B=(\Lambda(y),0)$ with $|y|=2k+1$, $k\ge1$, and
\[
H^*(P_gX;\Q)\cong\Lambda(y)\otimes\Q[\bar y],
\qquad \bar y=sy,\quad |\bar y|=2k.
\]

(iv) $\chi_o=\chi_e=1$. In this case
$B=(\Lambda(x,y),\delta)$ with $\delta x=0$, $\delta y=x^{n+1}$,
$|x|=d\in2\N$ and $|y|=d(n+1)-1$. Set $\bar x=sx$ and
$\bar y=sy$, the cohomology is
\bea
H^*(P_gX;\Q)\cong\Q\cdot1\oplus
\left(\frac{\Lambda^+(x,\bar x)}{(x^{n+1},x^n\bar x)}
\otimes\Q[\bar y]\right).
\eea
where $\Lambda^+(x,\bar x)$ is the augmentation ideal in the free graded-commutative algebra on $x, \bar{x}$, the positive-degree summand gives $\widetilde H^*(P_gX;\Q)$ which was used in Theorem 4.5 of \cite{GKV}.
}

\medskip

The transfer identifies the ordinary rational cohomology of $P_gX$
with that of $\Lambda_gM$. By the Classification Theorem and Proposition 3.2, the proof of Theorem 1.1
is divided into four cases.

{\bf Case 1.} {\it $\{b_q(\Lambda_gM)\}_{q\in\N}$ is unbounded.}

In this case, there are infinitely many non-contractible closed geodesics
on $(M,F)$, since otherwise the finiteness of the number of prime non-contractible closed geodesics, together with the mean index estimates, implies that the Morse numbers $M_q$ are uniformly bounded in sufficiently high degrees $q$,
note that $b_q(\Lambda_g M)\leq\bar \beta_q+\bar \beta_{q-1}\leq M_q+M_{q-1}$ by the Gysin sequence of the circle bundle:\[S^1\rightarrow ES^1\times\Lambda_g M\rightarrow ES^1\times_{S^1}\Lambda_g M \]
it's a contradiction to the unboundedness of the Betti numbers, see \cite{GKV} for the Riemannian metric case.

For the remaining cases, we argue by contradiction and assume that there is exactly one
non-contractible prime closed geodesic $c$ on $(M,F)$. Since $g^2=e$,
only the odd iterates $c^{2m-1}$ belong to $\Lambda_gM$. Moreover,
$E(c)$ is minimal on this component, and thus $i(c)=0$.
By Proposition 3.1 and Lemma 3.4, the Morse type numbers
\(
M_q=\sum_{m\ge1}\rank\overline C_q(E,c^{2m-1})
\)
satisfy $M_{2j-1}=0$ for every $j\in\N$, as stated in Lemma 4.1.
This gives a contradiction in Case (iv) of Proposition 3.2. In Cases (ii) and (iii) the odd
$S^1$-equivariant Betti numbers also vanish, and the Morse inequalities
give $M_q=\bar\beta_q$. In Section 4 we first finish Case (iv), then
consider Case (ii), and finally treat Case (iii) by systems of irrational numbers.

{\bf Case 2.} {\it $\chi_o=\chi_e=0$.}

In this case $B=\Q$. Proposition 3.2 gives
\begin{align}
\bar\beta_q=
\begin{cases}
1,&q\in2\N\cup\{0\},\\
0,&q\text{ odd},
\end{cases}\qquad
\bar B(\Lambda_gM;\Q)=1/2.
\label{1.8}
\end{align}
Under the assumption, the mean index equality in Theorem 3.7 gives $\hat i(c)=1$.

We apply the common index jump theorem (cf. Theorem 2.8) of Long and Zhu
\cite{LoZ} to this path $\gamma$ and its
second iteration $\gamma^2$, which yield
\bea
i(c^{2N-1})=2N-i(c),\qquad i(c^{2N+1})=2N+i(c).
\eea
Since $i(c)=0$, the two odd iterates have the same index $2N$. Each contributes one to $M_{2N}$,
while the Morse inequalities give
\(
1=\bar\beta_{2N}=M_{2N}\ge2.
\)
This contradiction proves Theorem 1.1 in this case.

{\bf Case 3.} {\it $\chi_o=1$ and $\chi_e=0$.}

In this case $B=(\Lambda(y),0)$ with $|y|=2k+1$, $k\in\N$.
Proposition 3.2 gives
\begin{equation}
\bar\beta_q=
\begin{cases}
2,&q\in\{2kj\mid j\in\N\},\\
1,&q\in(2\N\cup\{0\})\setminus\{2kj\mid j\in\N\},\\
0,&q\text{ odd}.
\end{cases}\qquad
\bar B(\Lambda_gM;\Q)=\frac{k+1}{2k}.
\label{1.10}
\end{equation}
Under the assumption, Theorem 3.7 gives $\hat i(c)=k/(k+1)$.
Using the precise index iteration formula and $i(c)=0$, the proof in
Section 4 obtains the following expressions in Lemma 4.2:
\bea
\begin{aligned}
i(c^{2m-1})
&=(2m-1)\frac{k}{k+1}+s
 -2\sum_{j=1}^{s}\{(2m-1)\hat\theta_j\},\\
\sum_{j=1}^{s}\hat\theta_j
&=\frac{s}{2}+\frac{k}{2(k+1)},
\end{aligned}
\eea
where $\hat\theta_j\notin\Q$, $s=2r\in2\N$ if $c$ is non-orientable, $s=r\in\N$ if $c$ is orientable. For
$m=2(k+1)l+2L+1$, formula (\ref{4.28}) describes $i(c^m)$ in terms
of the interval $I_i(L)$ containing
$\sum_{j=2}^{s}\{m\hat\theta_j\}$. The remaining argument follows
\cite{LiX,LLX} and is divided according to the rank of this system
of irrational numbers.

When $\rank\{\hat\theta_1,\ldots,\hat\theta_{s}\}=1$, set
$\hat\theta_j=p_j\theta+\xi_j$, where $p_j\in\Z\setminus\{0\}$,
$\theta\notin\Q$ and $\xi_j\in\Q$. We have
\bea
\sum_{j=1}^{s}p_j=0,\qquad
\left\{\sum_{j=1}^{s}\xi_j\right\}
=\{\frac{s}{2}\}+\frac{k}{2(k+1)}\in(0,1)\setminus\{1/2\}.
\eea
By Theorem 4.9, a rational translation of the system gives a nonzero
absolute difference number. The auxiliary functions $f_L$ in
(\ref{4.47}) and Lemma 4.10 then distinguish the index intervals for
$L=0$, while keeping the intervals the same for the relevant nonzero
values of $L$. Kronecker's approximation theorem supplies odd iterates
with parameters close to the two chosen endpoints. Lemma 4.3 restricts
the iterates contributing to a given index to finitely many offsets.
As in (\ref{4.54})--(\ref{4.60}), the contributions at these offsets
persist at the second parameter, which has one additional contribution
from its central iterate. The corresponding Betti numbers are equal
by their periodicity in positive even degrees. This contradicts
$M_q=\bar\beta_q$.

When the rank $\mathcal R$ is at least two, Lemma 4.11 provides a
change of variables for which every coefficient in one coordinate is
nonzero. We associate a rank-one system to this coordinate and apply
Theorem 4.9. The multivariable functions $g_L$ in (\ref{4.70}),
Kronecker's approximation theorem and Lemma 4.12 give the corresponding
interval comparison near $(0,\ldots,0)$ and $(1,0,\ldots,0)$.
The same counting argument yields a contradiction. This finishes
the proof in Case 3.

{\bf Case 4.} {\it $\chi_o=\chi_e=1$.}

In this case $B=(\Lambda(x,y),\delta)$ with
$\delta x=0$, $\delta y=x^{n+1}$ and $d=|x|\in2\N$.
It follows from Proposition 3.2 that $\bar\beta_{d-1}=1$. Since $d-1$ is odd,
Lemma 4.1 gives $M_{d-1}=0$, whereas the Morse inequality gives
$M_{d-1}\ge\bar\beta_{d-1}=1$. This contradiction proves Theorem 1.1 in this
case.

\medskip

In this paper, let $\N$, $\N_0$, $\Z$, $\Q$, and $\R$ denote the sets of natural integers,
non-negative integers, integers, rational numbers, and real numbers respectively. We use only singular
homology and co-homology modules with $\Q$-coefficients. We use also notations $[a]=\max\{k\in\Z\,|\,k\le a\}$,
$\mathcal{E}(a)=\min\{k\in\Z\,|\,k\ge a\}$, $\varphi(a)=\mathcal{E}(a)-[a]$, and $\{a\}=a-[a]$ for any
$a\in\R$.

\medskip

{\bf Acknowledgements.} The authors thank Hans-bert Rademacher and Iskander A. Taimanov for many valuable suggestions and comments for the original manuscript.

\medskip

{\bf AI disclosure.} The ideas of the proofs are based on our previous works. ChatGPT was only used to assist with literature searching and proof examining. The authors are responsible for all statements and proofs in this paper.

\setcounter{equation}{0}
\section{Index iteration theory for closed geodesics }%{Section 2}

In this section,  we review the index iteration theory developed by Y.  Long and his coworkers. In \cite{Lon1} of 1999, Y. Long established the basic normal form decomposition of symplectic matrices.
Based on this result he further established the precise iteration formulae of indices of symplectic
paths in \cite{Lon2} of 2000.

As usual, the symplectic group $\Sp(2n)$ is defined by
$$ \Sp(2n) = \{M\in {\rm GL}(2n,\R)\,|\,M^TJM=J\}, $$
whose topology is induced from that of $\R^{4n^2}$,
where  $J=\left(
              \begin{array}{cc}
                0 & -I_{n} \\
                I_{n} & 0 \\
              \end{array}
            \right)$ and
$I_n$ is the identity matrix in $\R^n$.
For $\tau>0$ we are interested in paths in $\Sp(2n)$:
$$ \P_{\tau}(2n) = \{\ga\in C([0,\tau],\Sp(2n))\,|\,\ga(0)=I_{2n}\}, $$
which is equipped with the topology induced from that of $\Sp(2n)$. The
following real function was introduced in \cite{Lon2}:
$$ D_{\om}(M) = (-1)^{n-1}\ol{\om}^n\det(M-\om I_{2n}), \qquad
          \forall \om\in\U,\, M\in\Sp(2n). $$
Thus for any $\om\in\U$ the following codimension $1$ hypersurface in $\Sp(2n)$ is
defined in \cite{Lon2}:
$$ \Sp(2n)_{\om}^0 = \{M\in\Sp(2n)\,|\, D_{\om}(M)=0\}.  $$
For any $M\in \Sp(2n)_{\om}^0$, we define a co-orientation of $\Sp(2n)_{\om}^0$
at $M$ by the positive direction $\frac{d}{dt}Me^{t\ep J}|_{t=0}$ of
the path $Me^{t\ep J}$ with $0\le t\le 1$ and $\ep>0$ being sufficiently
small. Let
\bea
\Sp(2n)_{\om}^{\ast} &=& \Sp(2n)\bs \Sp(2n)_{\om}^0,   \nn\\
\P_{\tau,\om}^{\ast}(2n) &=&
      \{\ga\in\P_{\tau}(2n)\,|\,\ga(\tau)\in\Sp(2n)_{\om}^{\ast}\}, \nn\\
\P_{\tau,\om}^0(2n) &=& \P_{\tau}(2n)\bs  \P_{\tau,\om}^{\ast}(2n).  \nn\eea
For any two continuous arcs $\xi$ and $\eta:[0,\tau]\to\Sp(2n)$ with
$\xi(\tau)=\eta(0)$, set $${\eta}*{\xi}(t)=\left\{ \begin{array}{ll}\xi(2t),~\mbox{if}~0\leq t\leq \tau/2,\\
\eta(2t-\tau),~\mbox{if}~\tau/2\leq t \leq \tau.\end{array}\right.$$ Given any
two $2m_k \times 2m_k$ square block matrices $M_{k}=\left(
              \begin{array}{cc}
                A_k & B_k \\
                C_k & D_k\\
              \end{array}
            \right)$
with $k = 1, 2$,
as in \cite{Lon3}, the $\;\diamond$-product of $M_1$ and $M_2$ is defined by
the following
$2(m_1+m_2)\times 2(m_1+m_2)$ matrix $M_1 \diamond M_2$:\\$$M_1 \diamond M_2=\left(
                            \begin{array}{cccc}
                              A_1& 0 & B_1 & 0 \\
                              0 & A_2 & 0 & B_2 \\
                              C_1 & 0 & D_1 & 0 \\
                              0 & C_2 & 0 & D_2 \\
                            \end{array}
                          \right)$$ %\diamond=\diamond
Denote by $M^{\diamond k}$ the $k$-fold $\diamond$-product $M\diamond\cdots\diamond M$. Note
that the $\diamond$-product of any two symplectic matrices is symplectic. For any two
paths $\ga_j\in\P_{\tau}(2n_j)$ with $j=0$ and $1$, let
$\ga_0\diamond\ga_1(t)= \ga_0(t)\diamond\ga_1(t)$ for all $t\in [0,\tau]$. Let

\be \xi_n(t)=\left(
                 \begin{array}{cc}
                   2-\frac{t}{\tau}& 0 \\
                   0 & (2-\frac{t}{\tau})^{-1} \\
                 \end{array}
               \right)^{\diamond n},~for~~0 \leq t \leq \tau.  \lb{2.1}\ee
{\bf Definition 2.1.} (cf. \cite{Lon3}) {\it For any $\om\in\U$ and
$M\in \Sp(2n)$, define
\be  \nu_{\om}(M)=\dim_{\C}\ker_{\C}(M - \om I_{2n}).  \lb{2.2}\ee
For any $\tau>0$ and $\ga\in \P_{\tau}(2n)$, define
\be  \nu_{\om}(\ga)= \nu_{\om}(\ga(\tau)).  \lb{2.3}\ee

If $\ga\in\P_{\tau,\om}^{\ast}(2n)$, define
\be i_{\om}(\ga) = [\Sp(2n)_{\om}^0: \ga\ast\xi_n],  \lb{2.4}\ee
where the right hand side of (\ref{2.4}) is the usual homotopy intersection
number, and the orientation of $\ga\ast\xi_n$ is its positive time direction under
homotopy with fixed end points.

If $\ga\in\P_{\tau,\om}^0(2n)$, we let $\mathcal{F}(\ga)$
be the set of all open neighborhoods of $\ga$ in $\P_{\tau}(2n)$, and define
\be i_{\om}(\ga) = \sup_{U\in\mathcal{F}(\ga)}\inf\{i_{\om}(\beta)\,|\,
                       \beta\in U\cap\P_{\tau,\om}^{\ast}(2n)\}.
               \lb{2.5}\ee
Then
$$ (i_{\om}(\ga), \nu_{\om}(\ga)) \in \Z\times \{0,1,\ldots,2n\}, $$
is called the index function of $\ga$ at $\om$. }

For any symplectic path $\ga\in\P_{\tau}(2n)$ and $m\in\N$,  we
define its $m$-th iteration $\ga^m:[0,m\tau]\to\Sp(2n)$ by
\be \ga^m(t) = \ga(t-j\tau)\ga(\tau)^j, \qquad
  {\rm for}\quad j\tau\leq t\leq (j+1)\tau,\;j=0,1,\ldots,m-1.
     \lb{2.6}\ee
We still denote the extended path on $[0,+\infty)$ by $\ga$.

{\bf Definition 2.2.} (cf. \cite{Lon3}) {\it For any $\ga\in\P_{\tau}(2n)$,
we define
\be (i(\ga,m), \nu(\ga,m)) = (i_1(\ga^m), \nu_1(\ga^m)), \qquad \forall m\in\N.
   \lb{2.7}\ee
The mean index $\hat{i}(\ga,m)$ per $m\tau$ for $m\in\N$ is defined by
\be \hat{i}(\ga,m) = \lim_{k\to +\infty}\frac{i(\ga,mk)}{k}. \lb{2.8}\ee
For any $M\in\Sp(2n)$ and $\om\in\U$, the {\it splitting numbers} $S_M^{\pm}(\om)$
of $M$ at $\om$ are defined by
\be S_M^{\pm}(\om)
     = \lim_{\ep\to 0^+}i_{\om\exp(\pm\sqrt{-1}\ep)}(\ga) - i_{\om}(\ga),
   \lb{2.9}\ee
for any path $\ga\in\P_{\tau}(2n)$ satisfying $\ga(\tau)=M$.}

For a given path $\gamma\in {\cal P}_{\tau}(2n)$ we consider to deform
it to a new path $\eta$ in ${\cal P}_{\tau}(2n)$ so that
\begin{equation}
i_1(\gamma^m)=i_1(\eta^m),\quad \nu_1(\gamma^m)=\nu_1(\eta^m), \quad
         \forall m\in {\bf N}, \label{2.10}
\end{equation}
and that $(i_1(\eta^m),\nu_1(\eta^m))$ is easy enough to compute. This
leads to finding homotopies $\delta:[0,1]\times[0,\tau]\to {\rm Sp}(2n)$
starting from $\gamma$ in ${\cal P}_{\tau}(2n)$ and keeping the end
points of the homotopy always stay in a certain suitably chosen maximal
subset of ${\rm Sp}(2n)$ so that (\ref{2.10}) always holds. This
set was first discovered in \cite{Lon2} as the path connected component
$\Omega^0(M)$ containing $M=\gamma(\tau)$ of the set
\begin{eqnarray}
  \Omega(M)=\{N\in{\rm Sp}(2n)\,&|&\,\sigma(N)\cap{\bf U}=\sigma(M)\cap{\bf U}\;
{\rm and}\;  \nonumber\\
 &&\qquad \nu_{\lambda}(N)=\nu_{\lambda}(M),\;\forall\,
\lambda\in\sigma(M)\cap{\bf U}\}. \label{2.11}
\end{eqnarray}
Here $\Omega^0(M)$ is called the {\it homotopy component} of $M$ in
${\rm Sp}(2n)$.

In \cite{Lon2}-\cite{Lon3}, the following symplectic matrices were introduced
as {\it basic normal forms}:
\begin{eqnarray}
D(\lambda)=\left(
                                                   \begin{array}{cc}
                                                      \lambda & 0 \\
                                                      0 & \lambda^{-1} \\
                                                    \end{array}
                                                  \right),
                                                  \lambda=\pm 2,\lb{2.12}\\
N_1(\lambda,b)=\left(
 \begin{array}{cc}
  \lambda & b \\
   0& \lambda \\
   \end{array}
   \right),\lambda=\pm1,b=\pm1,0, \lb{2.13}\\
R(\theta)=\left(
            \begin{array}{cc}
              \cos{\theta} & -\sin{\theta} \\
              \sin{\theta} & \cos{\theta} \\
            \end{array}
          \right),\theta\in (0,\pi)\cup(\pi,2\pi),
                     \lb{2.14}\\
N_2(\omega,B)=\left(
                \begin{array}{cc}
                  R(\theta) & b \\
                  0 &  R(\theta) \\
                \end{array}
              \right),\theta\in (0,\pi)\cup(\pi,2\pi),
                     \lb{2.15}\end{eqnarray}
where  $\omega=e^{\sqrt{-1}\theta}$, $b=\left(
           \begin{array}{cc}
             b_1 & b_2 \\
             b_3 & b_4 \\
           \end{array}
         \right)$ with  $b_i\in\R$ and  $b_2\not=b_3$. We call $ N_2(\omega, b)$ is
trivial if $(b_2-b_3)\sin\theta>0$;  it is
non-trivial if $(b_2-b_3)\sin\theta<0$.

Splitting numbers possess the following properties:

{\bf Lemma 2.3.} (cf. Lemma 9.1.5 of \cite{Lon3}) {\it Splitting
numbers $S_M^{\pm}(\om)$ are well defined, i.e., they are independent of the choice
of the path $\ga\in\P_\tau(2n)$ satisfying $\ga(\tau)=M$ appeared in (\ref{2.9}).
For $\om\in\U$ and $M\in\Sp(2n)$, splitting numbers $S_N^{\pm}(\om)$ are constant
for all $N\in\Om^0(M)$. }

{\bf Lemma 2.4.} (cf. List 9.1.12 of \cite{Lon3})
{\it For $M\in\Sp(2n)$ and $\om\in\U$, there hold
\begin{eqnarray}
S_M^{\pm}(\om) &=& 0, \qquad {\it if}\;\;\om\not\in\sg(M).  \lb{2.16}\\
S_{N_1(1,a)}^+(1) &=&\left\{
                \begin{array}{cc}
                1, &\quad {\rm if}\;\; a\ge 0, \cr
0, &\quad {\rm if}\;\; a< 0. \cr \end{array}\right. \lb{2.17}\eea

For any $M_i\in\Sp(2n_i)$ with $i=0$ and $1$, there holds }
\be S^{\pm}_{M_0\diamond M_1}(\om) = S^{\pm}_{M_0}(\om) + S^{\pm}_{M_1}(\om),
    \qquad \forall\;\om\in\U. \lb{2.18}\ee

We have the following

{\bf Theorem 2.5.} (cf. \cite{Lon2} and Theorem 1.8.10 of \cite{Lon3}) {\it For
any $M\in\Sp(2n)$, there is a path $f:[0,1]\to\Om^0(M)$ such that $f(0)=M$ and
\be f(1) = M_1\diamond\cdots\diamond M_k,  \lb{2.19}\ee
where each $M_i$ is a basic normal form listed in (\ref{2.12})-(\ref{2.15})
for $1\leq i\leq k$.}

The next theorem is the precise index iteration formulae  of Y. Long (cf. Theorem 8.3.1 and Corollary 8.3.2 of \cite{Lon3}).

{\bf Theorem 2.6.} {\it Let $\gamma\in\P_\tau(2n)$, Then there exists a path $f\in
C([0,1],\Omega^0(\gamma(\tau))$ such that $f(0)=\gamma(\tau)$ and
\bea f(1)=&&N_1(1,1)^{\diamond p_-} \diamond I_{2p_0}\diamond
N_1(1,-1)^{\diamond p_+}
\diamond N_1(-1,1)^{\diamond q_-} \diamond (-I_{2q_0})\diamond
N_1(-1,-1)^{\diamond q_+}\nn\\
&&\diamond R(\theta_1)\diamond\cdots \diamond R(\theta_{r'})\diamond  R(\theta_{r'+1})\diamond\cdots \diamond R(\theta_r)
\diamond N_2(\omega_1, u_1)\diamond\cdots\diamond N_2(\omega_{r_*}, u_{r_*}) \nn\\
&&\diamond N_2(\lm_1, v_1)\diamond\cdots\diamond N_2(\lm_{r_0}, v_{r_0})
\diamond M_0 \lb{2.20}\eea
where $ N_2(\omega_j, u_j) $s are
non-trivial and   $ N_2(\lm_j, v_j)$s  are trivial basic normal
forms; $\sigma (M_0)\cap U=\emptyset$; $p_-$, $p_0$, $p_+$, $q_-$,
$q_0$, $q_+$, $r$, $r_*$ and $r_0$ are non-negative integers;
$\omega_j=e^{\sqrt{-1}\alpha_j}$, $
\lambda_j=e^{\sqrt{-1}\beta_j}$; $\theta_j\in(\pi,2\pi)$ for $1\le j\le r^\prime$, $\theta_j\in(0,\pi)$ for $r^\prime< j\le r$, $\alpha_j$, $\beta_j$
$\in (0, \pi)\cup (\pi, 2\pi)$; The above parameters refer to the chosen normal form decomposition. Then
\bea i(\gamma, m)=&&m(i(\gamma,
1)+p_-+p_0-r)+2\sum_{j=1}^r\mathcal{E}\left(\frac{m\theta_j}{2\pi}\right)-r
-p_--p_0\nn\\&&-\frac{1+(-1)^m}{2}(q_0+q_+)+2\left(
\sum_{j=1}^{r_*}\varphi\left(\frac{m\alpha_j}{2\pi}\right)-r_*\right),
\lb{2.21}\eea
\bea \nu(\gamma, m)=&&\nu(\gamma,
1)+\frac{1+(-1)^m}{2}(q_-+2q_0+q_+)+2(r+r_*+r_0)\nn\\
&&-2\left(\sum_{j=1}^{r}\varphi\left(\frac{m\theta_j}{2\pi}\right)+
\sum_{j=1}^{r_*}\varphi\left(\frac{m\alpha_j}{2\pi}\right)
+\sum_{j=1}^{r_0}\varphi\left(\frac{m\beta_j}{2\pi}\right)\right),\lb{2.22}\eea
\bea \hat i(\gamma, 1)=i(\gamma, 1)+p_-+p_0-r+\sum_{j=1}^r
\frac{\theta_j}{\pi}.\lb{2.23}\eea
We have $i(\gamma, 1)$ is odd if $f(1)=N_1(1, 1)$, $I_2$,
$N_1(-1, 1)$, $-I_2$, $N_1(-1, -1)$ and $R(\theta)$; $i(\gamma, 1)$ is
even if $f(1)=N_1(1, -1)$ and $ N_2(\omega, b)$; $i(\gamma, 1)$ can be any
integer if $\sigma (f(1)) \cap \U=\emptyset$.}

We have the following properties in the index iteration theory.

{\bf Theorem 2.7. } (cf. Theorem 2.2 of \cite{Lon3})
{\it Let $\gamma\in\P_\tau(2n)$, then for any $m\in\N$, there holds
$$\nu(\gamma, m)-\frac{e(M)}{2}\le i(\gamma, m+1)-i(\gamma, m)-i(\gamma, 1)
\le \nu(\gamma, 1)-\nu(\gamma, m+1)+\frac{e(M)}{2}$$
where $M=\gamma(\tau)$ and $e(M)$ is the sum of the algebraic mutiplicities of the eigenvalues of $M$ on the unit circle. In particular, $0\le e(M)\le 2n$ for $M\in\mbox{Sp}(2n)$.}

The following is the common index jump theorem of Y. Long and C. Zhu.

{\bf Theorem 2.8.} (cf. Theorems 4.1-4.3 of \cite{LoZ}) {\it 	Let $\gamma_k\in\P_{\tau_k}(2n)$ for
$ k = 1,\ldots,q$ be a finite collection of symplectic paths. Let $M_k = \gamma_k(\tau_k)$.
Suppose $\hat i(\gamma_k, 1) > 0$,	for all $k = 1,\ldots,q$.
Then there exist infinitely many $(N,m_1,\ldots,m_q)\in\N^{q+1}$  such that
\bea \nu(\gamma_k, 2m_k -1)&=&	\nu(\gamma_k, 1),\lb{2.24}\\
\nu(\gamma_k, 2m_k +1)&=&	\nu(\gamma_k, 1),\lb{2.25}\\
i(\gamma_k, 2m_k -1)+\nu(\gamma_k, 2m_k -1)&=&	
2N-\left(i(\gamma_k, 1)+2S^+_{M_k}(1)-\nu(\gamma_k, 1)\right),\lb{2.26}\\
i(\gamma_k, 2m_k+1)&=&2N+i(\gamma_k, 1),\lb{2.27}\\
i(\gamma_k, 2m_k)&\ge&2N-\frac{e(M_k)}{2}\ge 2N-n,\lb{2.28}\\
i(\gamma_k, 2m_k)+\nu(\gamma_k, 2m_k)&\le&2N+\frac{e(M_k)}{2}\le 2N+n,\lb{2.29}
\eea
for every $k=1,\ldots,q$.
Moreover we have
\bea \min\left\{\left\{\frac{m_k\theta}{\pi}\right\},\,1-\left\{\frac{m_k\theta}{\pi}\right\}\right\}
<\delta,\lb{2.30}\eea
whenever $e^{\sqrt{-1}\theta}\in\sigma(M_k)$
and $\delta$ can be chosen as small as we want (cf. (4.43) of \cite{LoZ}). More precisely, by (4.10) and (4.40) in \cite{LoZ} , we have
\bea m_k=\left(\left[\frac{N}{M\hat i(\gamma_k, 1)}\right]+\xi_k\right)M,\quad 1\le k\le q,\lb{2.31}\eea
where $\xi_k=0$ or $1$ for $1\le k\le q$ and $\frac{M\theta}{\pi}\in\Z$
whenever $e^{\sqrt{-1}\theta}\in\sigma(M_k)$ and $\frac{\theta}{\pi}\in\Q$
for some $1\le k\le q$. By (4.20) in Theorem 4.1 of \cite{LoZ},
for any $\epsilon>0$, we can choose $N$ and $\{\xi_k\}_{1\le k\le q}$ such that
\bea \left|\frac{N}{M\hat i(\gamma_k, 1)}-\left[\frac{N}{M\hat i(\gamma_k, 1)}\right]-\xi_k\right|<\epsilon,\quad 1\le k\le q.\lb{2.32}\eea
Furthermore, given $M_0\in\N$, by the proof of Theorem 4.1 of \cite{LoZ} (Theorem 11.1.1 of \cite{Lon3}), we may  further require  $M_0|N$ (since the closure of the set $\{\{Nv\}: N\in\N, \;M_0|N\}$ is still a closed additive subgroup of $\T^h$ for some $h\in\N$, where we use notations  as (4.21)  in \cite{LoZ}.  Then we can use the proof of Step 2 in Theorem 4.1 of \cite{LoZ} to get $N$.)

In fact,  Let  $\mu_i=\sum_{\theta\in(0, 2\pi)}S_{M_i}^-(e^{\sqrt{-1}\theta})$ for $1\le i\le q$
and $\alpha_{i, j}=\frac{\theta_j}{\pi}$ where $e^{\sqrt{-1}\theta_j}\in\sigma(M_i)$
for $1\le j\le \mu_i$ and $1\le i\le q$.  Let $h=q+\sum_{1\le i\le q}\mu_i$ and
\bea  v=\left(\frac{1}{M\hat i(\gamma, 1)},\dots, \frac{1}{M\hat i(\gamma_q, 1)},
\frac{\alpha_{1, 1}}{\hat i(\gamma, 1)}, \frac{\alpha_{1, 2}}{\hat i(\gamma, 1)},
\dots\frac{\alpha_{1, \mu_1}}{\hat i(\gamma, 1)},
\frac{\alpha_{2, 1}}{\hat i(\gamma_2, 1)},\dots,
\frac{\alpha_{q, \mu_q}}{\hat i(\gamma_q, 1)}\right)\in\R^h.
\lb{2.33}\eea
Then  the above theorem is equivalent to find a vertex
\bea\chi=(\chi_1,\ldots,\chi_q,\chi_{1, 1},\chi_{1, 2},\ldots,\chi_{1, \mu_1},
\chi_{2, 1},\ldots,\chi_{q, \mu_q})\lb{2.34}\eea
of the cube $[0,\,1]^h$ and infinitely many integers $N\in\N$ such that
\bea|\{Nv\}-\chi|<\epsilon\lb{2.25}\eea
 for any given $\epsilon$ small enough (cf. P. 346 and 349 of \cite{LoZ}
).}

\medskip

From now on, we focus on the Morse index iteration formulae for non-orientable (or orientable) closed geodesics. The definition of Morse index $i(c)$ for a closed geodesic $c$ will be reviewed in the begining of Section 3.

\medskip

{\bf Lemma 2.9.} (cf. Theorem 1.1 of \cite{Liu1}) {\it  Let $c$ be a non-orientable closed geodesic on a $d$-dimensional Finsler manifold with its linear Poincar\'{e} map $P_c$.  Then, the following two claims hold.

$(i)$ If $d$ is even, there is a symplectic path $\gamma$ in Sp$(2d-2)$  with $\gamma(0)=I$ and $\gamma(1)=-P_c$ satisfying
$$  (i(c^m),\nu(c^{m}))= \left\{\begin{array}{ll}
    (i_{-1}(\gamma^m),\nu_{-1}(\gamma^{m})),\; &\ if\ m\ is\ odd,\\
    (i(\gamma^m),\nu(\gamma^{m}))\; &\ if \ m  \ is \ even.\\
    \end{array}\right. $$

$(ii)$ If $d$ is odd, there is a symplectic path $\td\gamma$ in Sp$(2d)$ with $\td\gamma(0)=I$ and $\td\gamma(1)=N_{1}(-1,1)\diamond (-P_c)$ satisfying
$$
    (i(c^m),\nu(c^{m}))= \left\{\begin{array}{ll}
    (i_{-1}(\td{\gamma}^m),\nu_{-1}(\td{\gamma}^{m})-1),\; &\ if\ m\ is\ odd,\\
    (i(\td{\gamma}^m),\nu(\td{\gamma}^{m})-1)\; &\ if \ m  \ is \ even.\\
    \end{array}\right. $$
}

{\bf Proof.}
Let $\mathcal N_c$ denote the normal bundle of $c$ with respect to the
fundamental tensor $g_{\dot c}$. Since $c$ is nonorientable and the tangent
line bundle $\mathbb R\dot c$ is trivial, both $\mathcal N_c$ and $c^*TM$
are nonorientable. Set
\[
\mathcal E=
\begin{cases}
\mathcal N_c,& d\text{ is even},\\
c^*TM,& d\text{ is odd},
\end{cases}
\qquad
\epsilon=
\begin{cases}
0,& d\text{ is even},\\
1,& d\text{ is odd}.
\end{cases}
\]
Then $\mathcal E$ has odd rank $r=d-1+\epsilon$. Every nonorientable real
vector bundle of odd rank over $S^1$ admits an orthonormal trivialization
over $\mathbb R$ whose transition matrix is $-I_r$. We therefore choose a
smooth $g_{\dot c}$-orthonormal frame
\begin{equation}\label{eq:anti-frame}
e(t)=(e_1(t),\ldots,e_r(t)),
\qquad e(t+1)=-e(t).
\end{equation}
In the full-bundle case, its initial value is chosen to be adapted to the
tangent--normal splitting.

Consider $c^m$ on the interval $[0,m]$. Writing a vector field as
$\xi(t)=e(t)u(t)$, its geometric periodicity becomes
\begin{equation}\label{eq:twisted-boundary}
u(m)=\omega_m u(0),
\qquad \omega_m=(-1)^m.
\end{equation}
In this frame, the index form restricted to $\mathcal E$ has the
Sturm--Liouville expression
\begin{equation}\label{eq:index-form}
\mathcal I_m(u,u)
=\int_0^m
\left(
\langle Pu',u'\rangle
+2\langle Qu,u'\rangle
+\langle Ru,u\rangle
\right)\,dt,
\end{equation}
where $P=P^T>0$, $R=R^T$, and all coefficient matrices are $1$-periodic.
The periodicity follows from \eqref{eq:anti-frame} and the periodicity of
the geometric Jacobi equation.

Put $p=Pu'+Qu$, and let $\Phi$ be the fundamental solution of the
associated Hamiltonian system in the phase coordinates $z=(p,u)$, with
$\Phi(0)=I_{2r}$. Since its coefficients are $1$-periodic, the fundamental
solution on $[0,m]$ is precisely the iterated path $\Phi^m$. For a Jacobi
field in the kernel of \eqref{eq:index-form}, the natural boundary
condition on the momentum is
\[
p(m)=\omega_m p(0).
\]
Thus the Jacobi boundary condition is $z(m)=\omega_m z(0)$.

The Morse--Maslov index theorem for positive Lagrangian systems (cf., Theorem 6.1 of \cite{Lon1}), with the
normalization of $i_\omega$ used in this paper, gives
\begin{equation}\label{eq:morse-maslov}
\operatorname{ind}\mathcal I_m=i_{\omega_m}(\Phi^m),
\qquad
\dim\ker\mathcal I_m=\nu_{\omega_m}(\Phi^m).
\end{equation}
The second equality also follows directly from the identification of
Jacobi fields with their initial Cauchy data.

The tangential and normal parts of the index form are orthogonal.
Moreover, the tangential form is a positive constant multiple of
$
\int_0^m |f'(t)|^2\,dt
$
, here \(f\) denotes the scalar coefficient of a tangential vector field \(\xi=f\dot c\), where \(f\in H^1([0,m],\mathbb R)\) satisfies \(f(m)=f(0)\).

It has index zero and a one-dimensional kernel.
Consequently,
\begin{equation}\label{eq:geometric-indices}
\operatorname{ind}\mathcal I_m=i(c^m),
\qquad
\dim\ker\mathcal I_m=\nu(c^m)+\epsilon.
\end{equation}
Rescaling $[0,m]$ to $[0,1]$ does not change these indices or nullities.

Suppose first that $d$ is even. Here $\mathcal E=\mathcal N_c$. Since
the terminal frame is the negative of the initial frame, the genuine
normal return map and the coordinate fundamental solution satisfy
\[
P_c=-\Phi(1).
\]
Taking $\gamma=\Phi$, we obtain $\gamma(1)=-P_c$.
Equations \eqref{eq:morse-maslov}--\eqref{eq:geometric-indices} yield
\[
(i(c^m),\nu(c^m))
=\bigl(i_{\omega_m}(\gamma^m),\nu_{\omega_m}(\gamma^m)\bigr),
\]
which proves (i).

Suppose now that $d$ is odd. In this case $\mathcal E=c^*TM$. The full
Jacobi return map splits into its tangential and normal parts. The
tangential Jacobi equation is $f''=0$, whose return matrix in the
coordinates $(p,q)=(f',f)$ is
$
U=\begin{pmatrix}1&0\\1&1\end{pmatrix}.
$
For
$
J_2=\begin{pmatrix}0&-1\\1&0\end{pmatrix},
$
we have
$
J_2UJ_2^{-1}
=\begin{pmatrix}1&-1\\0&1\end{pmatrix}
=N_1(1,-1).$

Hence the full geometric return map is symplectically conjugate to
$
N_1(1,-1)\diamond P_c.
$
By \eqref{eq:anti-frame}, its coordinate representation $\Phi(1)$ is
the negative of this return map. Therefore,
$
\Phi(1)\sim N_1(-1,1)\diamond(-P_c),
$
where $\sim$ denotes symplectic conjugacy.

After applying the corresponding constant symplectic conjugation to
$\Phi$, we obtain a path $\widetilde\gamma$ satisfying
\[
\widetilde\gamma(0)=I_{2d},
\qquad
\widetilde\gamma(1)=N_1(-1,1)\diamond(-P_c).
\]
Constant symplectic conjugation preserves the Maslov-type indices and
nullities, including those of all iterates. Thus
\eqref{eq:morse-maslov}--\eqref{eq:geometric-indices} imply
\[
(i(c^m),\nu(c^m))
=\bigl(
i_{\omega_m}(\widetilde\gamma^m),
\nu_{\omega_m}(\widetilde\gamma^m)-1
\bigr).
\]
Taking $\omega_m=-1$ for odd $m$ and $\omega_m=1$ for even $m$ proves (ii).

For any $m\in\N$, we define
$$E_{m}(a)=E\left(a-{1-(-1)^{m}\over4}\right),\ \varphi_{m}(a)=\varphi\left(a-{1-(-1)^{m}\over4}\right),\ \forall a\in\R.$$

The precise iteration formulae of Morse indices
for a non-orientable (or orientable) closed geodesic are the following.

{\bf Lemma 2.10.} (cf. Theorem 3.1 of \cite{LiX}) {\it Let $c$ be a non-orientable closed geodesic on a $d$-dimensional Finsler Manifold $M$ with its linear
Poincar$\acute{e}$ map $P_c$ and $-P_c$ can be continuously connected in $\Omega^0(-P_c)$ to the form (\ref{2.20}). Then for every $m\in2\N-1$, we have
\bea
i(c^{m})&=&m(i(c)+q_{0}+q_{+}-2r')-(q_{0}+q_{+})\nn\\
& &+2\sum_{j=1}^{r} E_{m}\left({m\theta_{j}\over2\pi}\right)+2\sum_{j=1}^{r_*} \varphi_{m}\left({m\alpha_{j}\over2\pi}\right)-2r_*,
\lb{ind1}\eea
and
\bea
\nu(c^{m})&=&\nu(c)+2\td{\varsigma}(c,m),\lb{null1}
\eea
where we denote by
$$\td{\varsigma}(c,m)=\left(r-\sum_{j=1}^{r}\varphi_{m}\left({m\theta_{j}\over2\pi}\right)\right)
+\left(r_{*}-\sum_{j=1}^{r_{*}}\varphi_{m}\left({m\alpha_{j}\over2\pi}\right)\right)
+\left(r_{0}-\sum_{j=1}^{r_{0}}\varphi_{m}\left({m\beta_{j}\over2\pi}\right)\right).$$
}

{\bf Proof.}  We only prove the case when $d$ is even and $m$ is odd, the other case can be proved similarly. By Lemma 2.9,  there exists a symplectic path $\gamma$ in Sp$(2d-2)$ with $\ga(0)=I$ and $\gamma(1)=-P_c$ such that
\bea (\text{ind}(c^{m}),\text{null}(c^{m}))=(i_{-1}(\gamma^{m}),\nu_{-1}(\gamma^{m})), \forall\ m\in2\N-1.
\lb{201679a}\eea
It together with the Bott-type formulae (cf. Theorem 9.2.1 of \cite{Lon3}) and Theorem 2.6 gives
\bea i_{-1}(\gamma^{m})&=&i(\gamma^{2m})-i(\gamma^{m})\nn\\
&=&m(i(\ga)+p_-+p_0-r )-(q_0+q_+)    \nn\\
&&+ 2\sum_{j=1}^r\left[{E}\left(\frac{m\th_j}{\pi}\right)-{E}\left(\frac{m\th_j}{2\pi}\right)\right]+ 2\sum_{j=1}^{r_{\ast}}\left[\vf\left(\frac{m\aa_j}{\pi}\right)-\vf\left(\frac{m\aa_j}{2\pi}\right)\right]\nn\\
&=&m(i(\ga)+p_-+p_0-r )-(q_0+q_+)  \nn\\
&&+ 2\sum_{j=1}^r{E}\left(\frac{m\th_j}{2\pi}-{1\over2}\right)+ 2\sum_{j=1}^{r_{\ast}} \vf\left(\frac{m\aa_j}{2\pi}-{1\over2}\right)-2r_{*}\lb{201679b}\\
&=&m(i_{-1}(\gamma)+q_0+q_+-2r' )-(q_0+q_+) \nn\\
&&+ 2\sum_{j=1}^r{E}\left(\frac{m\th_j}{2\pi}-{1\over2}\right)+ 2\sum_{j=1}^{r_{\ast}} \vf\left(\frac{m\aa_j}{2\pi}-{1\over2}\right)-2r_{*},\nn
\eea
where the third identity we have used $E(2a)-E(a)=E\left(a-{1\over2}\right)$ and $\vf(2a)-\vf(a)=\vf\left(a-{1\over2}\right)-1,$ and the last identity is due to
$$i(\gamma)=i_{-1}(\gamma)+(q_{0}+q_{+})+(r-2r')-(p_{0}+p_{-}),$$ which is a result of direct computation on splitting numbers based on Theorem 12.2.3 of \cite{Lon3}.

Observing by definition $\nu_{-1}(\gamma)=q_{-}+2q_{0}+q_{+}$, by Theorem 2.6 we obtain
 \bea \nu_{-1}(\gamma^{m})&=&\nu(\gamma^{2m})-\nu(\gamma^{m})\nn\\
&=&(q_{-}+2q_{0}+q_{+})-2\sum_{j=1}^r\left[{\varphi}\left(\frac{m\th_j}{\pi}\right)-{\varphi}\left(\frac{m\th_j}{2\pi}\right)\right]\nn\\
&&-2\sum_{j=1}^{r_{\ast}}\left[\vf\left(\frac{m\aa_j}{\pi}\right)-\vf\left(\frac{m\aa_j}{2\pi}\right)\right]
-2\sum_{j=1}^{r_{0}}\left[\vf\left(\frac{m\bb_j}{\pi}\right)-\vf\left(\frac{m\bb_j}{2\pi}\right)\right]\nn\\
&=&\nu_{-1}(\gamma)+2\left(r-\sum_{j=1}^r{\varphi}\left(\frac{m\th_j}{2\pi}-{1\over2}\right)\right)\lb{201679c} \\
&& +2\left(r_{*}-\sum_{j=1}^{r_{\ast}}\vf\left(\frac{m\aa_j}{2\pi}-{1\over2}\right)\right)
+2\left(r_0-\sum_{j=1}^{r_{0}}\vf\left(\frac{m\bb_j}{2\pi}-{1\over2}\right)\right).\nn
\eea
Thus (\ref{ind1}) and (\ref{null1})  immediately follow from (\ref{201679a}), (\ref{201679b}) and (\ref{201679c}). \hfill$\Box$

\medskip

{\bf Lemma 2.11.} (cf. Lemma 3.1 of \cite{LiX}) {\it Let $c$ be an orientable closed geodesic on an $d$-dimensional Finsler manifold with its Poincar\'{e} map $P_c$.  Then, there
exists a continuous symplecitic path  $\ga$ with $\ga(0)=I$ and $\ga(1)=P_c$ which can be continuously connected in $\Omega^0(P_c)$ to the form (\ref{2.20}) such that
\bea i(c^{m})=i(\ga^m)
&=& m(i(\ga)+p_-+p_0-r )  -(p_- + p_0+r) - {{1+(-1)^m}\over 2}(q_0+q_+) \nn\\&&
              + 2\sum_{j=1}^r{E}\left(\frac{m\th_j}{2\pi}\right)+ 2\sum_{j=1}^{r_{\ast}}\vf\left(\frac{m\aa_j}{2\pi}\right) - 2r_{\ast},
\lb{preind}\eea
and
\bea
{\rm null}(c^{m})=\nu(\ga^m)
 &=& \nu(\ga) + {{1+(-1)^m}\over 2}(q_-+2q_0+q_+) + 2\vs(c,m),    \lb{pnul}\eea
where we denote by}
\be \vs(c,m) = \left(r - \sum_{j=1}^r\vf\left(\frac{m\th_j}{2\pi}\right)\right)
             + \left(r_{\ast} - \sum_{j=1}^{r_{\ast}}\vf\left(\frac{m\aa_j}{2\pi}\right)\right)
             + \left(r_0 - \sum_{j=1}^{r_0}\vf\left(\frac{m\bb_j}{2\pi}\right)\right).    \nn\lb{yuxiang}\ee

\setcounter{equation}{0}
\section{Morse theory for non-contractible closed geodesics}%{Section 3}

In this section, we review the variational structures of closed geodesic,
all the details can be found in \cite{Rad2} or \cite{BaL}.

On a compact Finsler manifold $(M,F)$, we choose an auxiliary Riemannian
metric. This endows the space $\Lambda=\Lambda M$ of $H^1$-maps
$\gamma:S^1\rightarrow M$ with a natural structure of Riemannian
Hilbert manifold on which the group $S^1=\R/\Z$ acts continuously
by isometries, cf. \cite{Kli2}, Chapters 1 and 2. This action is
defined by translating the parameter, i.e.
$$ (s\cdot\gamma)(t)=\gamma(t+s) \qquad $$
for all $\gamma\in\Lm$ and $s,t\in S^1$.
The Finsler metric $F$ defines an energy functional $E$ and a length
functional $L$ on $\Lambda$ by
\be E(\gamma)=\frac{1}{2}\int_{S^1}F(\dot{\gamma}(t))^2dt,
 \quad L(\gamma) = \int_{S^1}F(\dot{\gamma}(t))dt.  \lb{3.1}\ee
Both functionals are invariant under the $S^1$-action. The critical points
of $E$ of positive energies are precisely the closed geodesics $c:S^1\to M$
of the Finsler structure. If $c\in\Lambda$ is a closed geodesic, then $c$ is
a regular curve, i.e. $\dot{c}(t)\not= 0$ for all $t\in S^1$, and this implies
that the second differential $E''(c)$ of $E$ at $c$ exists.
As usual we define the index $i(c)$ of $c$ as the maximal dimension of
subspaces of $T_c \Lambda$ on which $E^{\prime\prime}(c)$ is negative definite, and the
nullity $\nu(c)$ of $c$ so that $\nu(c)+1$ is the dimension of the null
space of $E^{\prime\prime}(c)$.

For $m\in\N$ we denote the $m$-fold iteration map
$\phi^m:\Lambda\rightarrow\Lambda$ by \be \phi^m(\ga)(t)=\ga(mt)
\qquad \forall\,\ga\in\Lm, t\in S^1. \lb{3.2}\ee We also use the
notation $\phi^m(\gamma)=\gamma^m$. For a closed geodesic $c$, the
mean index is defined to be: \be
\hat{i}(c)=\lim_{m\rightarrow\infty}\frac{i(c^m)}{m}. \lb{3.3}\ee

If $\gamma\in\Lambda$ is not constant then the multiplicity
$m(\gamma)$ of $\gamma$ is the order of the isotropy group $\{s\in
S^1\mid s\cdot\gamma=\gamma\}$. If $m(\gamma)=1$ then $\gamma$ is
called {\it prime}. Hence $m(\gamma)=m$ if and only if there exists a
prime curve $\tilde{\gamma}\in\Lambda$ such that
$\gamma=\tilde{\gamma}^m$.

For a closed geodesic $c$ we set $ \Lm(c)=\{\ga\in\Lm\mid E(\ga)<E(c)\}. $
If $A\subseteq\Lm$ is invariant under some subgroup $\Gamma$ of $S^1$,
we denote by $A/\Gamma$ the
quotient space of $A$ with respect to the action of $\Gamma$.

We will use the following critical $\Q$-module of a closed geodesic
$c\in\Lambda$:
\be \overline{C}_*(E,c)
   = H_*\left((\Lm(c)\cup S^1\cdot c)/S^1,\Lm(c)/S^1\right). \lb{3.4}\ee

We call a closed geodesic satisfying the isolation condition, if
the following holds:

{\bf (Iso)  For all $m\in\N$ the orbit $S^1\cdot c^m$ is an
isolated critical orbit of $E$. }

Note that if the number of prime closed geodesics on a Finsler manifold
is finite, then all the closed geodesics satisfy (Iso).

The following Propositions were proved in \cite{Rad2} and \cite{BaL}.

{\bf Proposition 3.1.} (cf. Satz 6.11 of \cite{Rad2} or Proposition
3.12 of \cite{BaL}) {\it Let $c$ be a prime closed geodesic on a
bumpy Finsler manifold $(M,F)$ satisfying (Iso). Then we have
\bea \overline{C}_q( E,c^m) = \left\{\begin{array}{ll}
     \Q, &\quad {\it if}\;\; i(c^m)-i(c)\in 2\Z,\;{\it and}\;
                   q=i(c^m)\;  \\
     0, &\quad {\it otherwise}. \end{array}\right.  \lb{3.5}\eea}

{\bf Proposition 3.2.}  {\it   Let $(M, F)$ be a compact bumpy Finsler manifold with $\pi_1(M)=\Z_2$, we have the following classification where we use the notations as in Classification Theorem of Section 1.

(i)   $\{b_i(\Lambda_g M)\}_{i\in\N}$ is unbounded.

(ii)  $\chi_o=\chi_e=0$,  thus $\Lambda Y=\Q$ and $H^\ast(\Lambda_g M)=\Q$.

In this case, the $S^1$-equivariant Poincar$\acute{e}$ series
of $\Lambda_g M$ is given by
\bea P^{S^1}(\Lambda_g M; \Q)(t)=\frac{1}{1-t^2}=1+t^2+t^4+\cdots \lb{3.6}\eea
Thus the Betti numbers
\bea \bar{\beta}_q=\rank H_q^{S^1}(\Lambda_g M;\Q)=\left\{\begin{array}{ll}
1,&\ if\ q\in 2\N\cup\{0\},\\
0,&\ otherwise.\\
\end{array}
\right.\lb{3.7}\eea
and the average $S^1$-equivariant Betti number of $\Lambda_g M$ satisfies
\be  \bar{B}(\Lambda_g M;\Q)\equiv\lim_{q\to+\infty}\frac{1}{q}\sum_{k=0}^{q}(-1)^k\bar{\beta}_k= \frac{1}{2}.   \lb{3.8}\ee

(iii) $\chi_o=1$, $\chi_e=0$, thus $\Lambda Y=\Lambda(y)$ and $H^\ast(\Lambda_g M)=\Lambda(y)\otimes\Q[\bar y]$.

In this case,  the $S^1$-equivariant Poincar$\acute{e}$ series
of $\Lambda_g M$ is given by
\bea P^{S^1}(\Lambda_g M; \Q)(t)&=&\frac{1-t^{2k+2}}{(1-t^2)(1-t^{2k})}=\frac{1}{1-t^2}+\frac{t^{2k}}{1-t^{2k}}\nn\\
&=&(1+t^2+t^4+\cdots+t^{2k}+\cdots)+(t^{2k}+t^{4k}+t^{6k}+\cdots).\lb{3.9}\eea
Thus the Betti numbers
\bea \bar{\beta}_q=\rank H_q^{S^1}(\Lambda_g M;\Q)=\left\{\begin{array}{ll}
    2,&\ if\ q\in \{(2k)j\mid j\in\N\},\\
    1,&\ if\ q\in (2\N\cup\{0\})\backslash\{(2k)j\mid j\in\N\},\\
    0,&\ otherwise.\\
    \end{array}\right. \lb{3.10}\eea
and the average $S^1$-equivariant Betti number of $\Lambda_g M$ satisfies
\be  \bar{B}(\Lambda_g M;\Q)\equiv\lim_{q\to+\infty}\frac{1}{q}\sum_{k=0}^{q}(-1)^k\bar{\beta}_k= \frac{k+1}{2k}.   \lb{3.11}\ee

(iv) $\chi_o=\chi_e=1$. In this case $B=(\Lambda(x,y),\delta)$ with
$\delta x=0$ and $\delta y=x^{n+1}$, and
\[
H^*(\Lambda_gM;\Q)\cong\Q\cdot1\oplus
\left(\frac{\Lambda^+(x,\bar x)}{(x^{n+1},x^n\bar x)}
\otimes\Q[\bar y]\right).
\]
where $\Lambda^+(x,\bar x)$  is the augmentation ideal in the free graded-commutative algebra on $x, \bar{x}$, the positive-degree summand is $\widetilde H^*(\Lambda_gM;\Q)$ which was used in Theorem 4.5 of \cite{GKV}.

In this case, the $S^1$-equivariant Poincar$\acute{e}$ series
of $\Lambda_g M$ is given by
\bea P^{S^1}(\Lambda_g M; \Q)(t)&=&\frac{1}{1-t^2}+\frac{t^{d-1}}{1-t^{d(n+1)-2}}\frac{1-t^{dn}}{1-t^d}, \lb{3.12}\eea
where $d=\deg x\in 2\N$.}

{\bf Proof.} We use continuous loop spaces, without changing the rational
Borel cohomology of the Sobolev loop spaces in the statement.
Put $X=\tilde M$ and
$$
P_gX=\{\gamma:\R\to X\mid\gamma(t+1)=g\gamma(t)\}.
$$

\smallskip
The covering projection $X\to M$ induces the following covering
\bea
p_*:P_gX\longrightarrow\Lambda_gM,
\qquad
\Lambda_gM\cong P_gX/\langle g\rangle.
\eea
Its deck transformation is $\gamma\mapsto g\gamma=\gamma(\,\cdot+1)$,
which is homotopic to the identity through parameter translations.
Therefore, we get the following
\bea
p_*^*:H^*(\Lambda_gM;\Q)\xrightarrow{\cong}H^*(P_gX;\Q).
\eea

Let $T=\R/2\Z$ and
$\rho:T\to S^1=\R/\Z$, $[s]_2\mapsto[s]_1$.
Parameter translation defines a $T$-action on $P_gX$, and $p_*$ is
$\rho$-equivariant. Since $(B\rho)^*: H^*(BS^1; \Q)\rightarrow H^*(BT; \Q) $ is multiplication by $2$ in degree $2$, where $BS^1=ES^1/S^1\simeq \mathbb{C}P^\infty$, $ES^1$ is contractible, and
\[
H^*(BS^1;\Q)=\Q[u],\qquad |u|=2.
\]

Note that $(B\rho)^*$ is an isomorphism of cohomology rings. We use the following Borel comparison principle: Suppose $f:U\to V$ is equivariant with respect to a compact-group homomorphism $\rho:T\to S^1$, the spaces $BT,BS^1$ are simply connected, and both $f$ and $B\rho$ induce rational cohomology isomorphisms. Then $ U_{hT}\equiv ET\times_TU\longrightarrow V_{hS^1}\equiv ES^1\times_{S^1}V$ induces a rational cohomology
isomorphism in the contravariant direction
\bea
H^*_{S^1}(V;\Q)\xrightarrow{\cong}H^*_T(U;\Q),
\eea
where
\bea
H^*_{S^1}(V;\Q)=H^*(V_{hS^1};\Q),\qquad H^*_T(U;\Q)=H^*(U_{hT};\Q).
\eea
In particular, for $U=P_gX$ and $V=\Lambda_gM$, this gives
$H^*_{S^1}(\Lambda_gM;\Q)\cong H^*_T(P_gX;\Q)$.

Choose a $g$-equivariant minimal Sullivan model (cf. \cite{FHT})
\bea
\varphi:A=(\Lambda V,d)\xrightarrow{\simeq}A_{PL}(X),\qquad g_A^2=\id,\qquad g_A(V)=V.
\eea
Such a model is obtained by averaging the linear splittings in the
inductive minimal model construction. Write $V=V^+\oplus V^-$ for
the eigenspace decomposition of $g_A$. Since $d(V^-)$ is anti-invariant, the ideal $I=(V^-)$ is
$d$-stable. Hence we get
\bea
B=A/I=(\Lambda Y,\delta),\qquad Y=V^+,
\qquad q:A\longrightarrow B,\qquad qg_A=q.
\eea

Let $X_A=|\langle A\rangle|$, $Z=|\langle B\rangle|$, and
$j=|q|:Z\to X_A$. The Sullivan realization gives $g_Aj=j$ (cf. \cite{Sullivan}).
Using $X_c=|\operatorname{Sing}X|$ and the map
adjoint to $\varphi$, we obtain a strictly equivariant zigzag
$$
X\xleftarrow{\varepsilon}X_c\xrightarrow{r}X_A
\xleftarrow{j}Z,
$$
where $\varepsilon$ is a weak equivalence, $r$ is a rational
equivalence, and the involution on $Z$ is trivial. It induces
$T$-equivariant maps
\bea
P_gX\xleftarrow{P\varepsilon}P_gX_c
\xrightarrow{Pr}P_{g_A}X_A
\xleftarrow{F_j}\Lambda Z,
\qquad F_j(\gamma)=j\circ\gamma,
\eea
where the $T$-action on $\Lambda Z$ factors through $\rho$.

We verify that these maps induce rational cohomology isomorphisms.
For a simply connected space $W$ and an involution $h$, the evaluation
fibration $\Omega W\to P_hW\to W$ has connecting homomorphism
$h_*-\id$ (cf. Theorem A of \cite{Taimanov}).
With $W_q=\pi_q(W)\otimes\Q$, its exact sequence gives, for $q\ge2$,
\bea
\begin{aligned}
0\longrightarrow
\mbox{coker}(h_*-\id:W_{q+1}\to W_{q+1})
&\longrightarrow\pi_q(P_hW)\otimes\Q\\
&\longrightarrow\ker(h_*-\id:W_q\to W_q)
\longrightarrow0,
\end{aligned}
\eea
and
\bea
\pi_1(P_hW)\otimes\Q
\cong \mbox{coker}(h_*-\id:W_2\to W_2).
\eea
The group $\pi_1(P_hW)$ is a quotient of the abelian group
$\pi_2(W)$. Moreover, its action on $\pi_q(P_hW)$ has the filtration
$$
0\subset
\im\bigl(\pi_q(\Omega W)\to\pi_q(P_hW)\bigr)
\subset\pi_q(P_hW),
$$
with trivial action on both successive quotients: the first assertion
follows from the simplicity of the connected $H$-space $\Omega W$
and the surjection $\pi_1(\Omega W)\to\pi_1(P_hW)$; the second
follows from $\pi_1(W)=0$. Thus $P_hW$ is nilpotent.

By construction,
\bea
j_*:\pi_q(Z)\xrightarrow{\cong}
\bigl(\pi_q(X_A)\bigr)^{g_A}.
\eea
For an involution on a rational vector space, the inclusion of the
fixed subspace induces an isomorphism onto the cokernel of
$h_*-\id$ as well. The evaluation exact sequences therefore show
that $F_j$ induces isomorphisms on all rational homotopy groups.
The same argument applies to $P\varepsilon$ and $Pr$.
Nilpotence and the rational Whitehead theorem \cite{FHT} now give
rational cohomology isomorphisms throughout the zigzag.
Comparing the Borel spectral sequences over $BT$, and then using
$\rho$, we have
\bea
H^*(\Lambda_gM; \Q)\cong H^*(\Lambda Z; \Q),
\qquad
H^*_{S^1}(\Lambda_gM; \Q)\cong H^*_{S^1}(\Lambda Z; \Q).
\eea

\smallskip
By the boundedness criterion and the classification of
Grove--Halperin--Vigu\'e-Poirrier (cf. Theorem 4.5 of \cite{GKV}),
$$
\{b_q(P_gX)\}_{q\ge0}\text{ is bounded}
\quad\Longleftrightarrow\quad
\chi_e\le\chi_o\le1,
$$
where $\chi_e=\dim Y^{\mathrm{even}}$ and
$\chi_o=\dim Y^{\mathrm{odd}}$.
An auxiliary Riemannian metric on $X$ may be averaged over the deck
group to make $g$ an isometry. This verifies the hypothesis used for
the topological classification in \cite{GKV}. The transfer identifies
the ordinary rational Betti numbers of $P_gX$ and $\Lambda_gM$,
and hence identifies their boundedness conditions. The averaging is
used only for this topological classification; the prescribed
Finsler metric $F$ is kept fixed in the variational problem considered
in Section~4. In the bounded cases the
possible models are
$$
\begin{array}{ll}
\mathrm{(ii)}&B=\Q;\\[2pt]
\mathrm{(iii)}&B=(\Lambda(y),0),\quad |y|=2k+1,\quad k\ge1;\\[2pt]
\mathrm{(iv)}&B=(\Lambda(x,y),\delta),\quad
\delta x=0,\quad\delta y=x^{n+1},\\
&|x|=d\in2\N,\quad |y|=d(n+1)-1,\quad n\ge1.
\end{array}
$$

The free-loop and equivariant model theorems
\cite{VPS,VPB} identify $H^*_{S^1}(\Lambda Z)$ with the cohomology of
\bea
\mathcal C_B=
\bigl(\Q[u]\otimes\Lambda(Y\oplus sY),D\bigr),
\qquad |u|=2,\quad |sv|=|v|-1,
\eea
where $s$ is the degree $-1$ derivation with
$s(v)=sv$, $s(sv)=0$, and
\bea
Du=0,\qquad Dv=\delta v+u\,sv,
\qquad D(sv)=-s(\delta v).
\eea
Setting $u=0$ gives the ordinary free-loop model.

\smallskip
\textbf{Proof of Case (ii).}
Here $H^*(\Lambda_gM; \Q)=\Q$ and
$\mathcal C_B=(\Q[u],0)$. Hence
\bea
P^{S^1}(\Lambda_gM;\Q)(t)=\frac1{1-t^2},
\qquad
\bar\beta_q=
\begin{cases}
1,&q\in2\N_0,\\
0,&q\text{ odd},
\end{cases}
\eea
and
\bea
\bar B(\Lambda_gM;\Q)
=\lim_{N\to\infty}\frac{[N/2]+1}{N}
=\frac12.
\eea
This proves (\ref{3.6})-(\ref{3.8}).

\smallskip
\textbf{Proof of Case (iii).}
Put $z=sy$, so $|z|=2k$. The model has zero differential,
and therefore
\bea
H^*(\Lambda_gM)\cong\Lambda(y)\otimes\Q[z].
\eea
The equivariant model is
\bea
\mathcal C_B=\Q[u,z]\otimes\Lambda(y),
\qquad Du=Dz=0,\qquad Dy=uz.
\eea
Since multiplication by $uz$ is injective on $\Q[u,z]$, the displayed
complex has $H^*(\mathcal C_B)\cong\Q[u,z]/(uz)$. For the argument
below, we use its graded vector-space decomposition
\bea
H^*(\mathcal C_B)\cong
\Q\cdot1\oplus\bigoplus_{j\ge1}\Q u^j
\oplus\bigoplus_{j\ge1}\Q z^j,
\qquad |u|=2,\quad |z|=2k.
\eea

Thus
\bea
P^{S^1}(\Lambda_gM;\Q)(t)
=\frac1{1-t^2}+\frac{t^{2k}}{1-t^{2k}}
=\frac{1-t^{2k+2}}{(1-t^2)(1-t^{2k})},
\eea
\bea
\bar\beta_q=
\begin{cases}
2,&q\in\{2kj:j\in\N\},\\
1,&q\in2\N_0\setminus\{2kj:j\in\N\},\\
0,&q\text{ odd},
\end{cases}
\eea
and
\bea
\begin{aligned}
\bar B(\Lambda_gM;\Q)
=\lim_{N\to\infty}
\frac{[N/2]+1+[N/(2k)]}{N}
=\frac12+\frac1{2k}=\frac{k+1}{2k}.
\end{aligned}
\eea
This proves (\ref{3.9})-(\ref{3.11}).

\smallskip
\textbf{Proof of Case (iv).}
Set $a=sx$, $b=sy$, and $\ell=d(n+1)-2$. Then
$|a|=d-1$, $|b|=\ell$, and the equivariant model becomes
\bea
&&\mathcal C=\Q[u,x,b]\otimes\Lambda(a,y),\\
&&Du=Da=0,\quad Dx=ua,\quad
Dy=x^{n+1}+ub,\quad Db=-(n+1)x^na.
\eea
Let $F=x^{n+1}+ub$ and $R=\Q[u,x,b]/(F)$.
Since $DF=0$, the quotient projection
\bea \mathcal C\longrightarrow\mathcal C'=R\otimes\Lambda(a)\eea
is a chain map. We claim that it is a quasi-isomorphism. In fact, give both complexes the finite decreasing filtration
\bea
F^0\mathcal C=\mathcal C,\quad F^1\mathcal C=a\mathcal C,\quad F^2\mathcal C=0,
\eea
and the analogous filtration on the quotient. The differential preserves it. In the associated graded complex, $Dx$ and $Db$ increase the number of $a$-factors and therefore do not appear in the filtration-zero differential. The only remaining nonzero generator differential is $Dy=F$. If $P=\Q[u,x,b]$, the $y$-part of the filtration-zero complex is the two-term complex $Py\xrightarrow{\cdot F}P$. Since $P$ is an integral domain and $F\ne0$, then the multiplication by $F$ is injective. Its cohomology is $P/(F)$, with no residual $y$-cohomology. Including the exterior factor $a$, the $E_1$ page of the original complex is $R\otimes\Lambda(a)$,
which is exactly the $E_1$ page of the quotient complex, and the projection induces this isomorphism. Since the filtration has only two steps, then the quotient projection is a quasi-isomorphism.

Define
\bea
L\equiv u\frac{\partial}{\partial x}
-(n+1)x^n\frac{\partial}{\partial b}.
\eea
As $L(F)=0$, this derivation descends to $R$, and the differential
on $\mathcal C'$ satisfies $Df=aL(f)$ and $D(af)=0$ for $f\in R$.
All elements of $R$ have even degree. Hence
\bea
H^{\mathrm{even}}(\mathcal C')=\ker L,
\quad
H^{\mathrm{odd}}(\mathcal C')=a\bigl(R/L(R)\bigr).
\eea
Because $F$ is monic in $x$, a $\Q$-basis of $R$ is
\bea
u^r x^i b^j,
\qquad r,j\ge0,\quad0\le i\le n.
\eea
Using $x^{n+1}=-ub$, we obtain
\bea
\begin{aligned}
L(u^r x^i b^j)
&=\bigl(i+(n+1)j\bigr)u^{r+1}x^{i-1}b^j,
&&1\le i\le n,\\
L(u^r b^j)
&=-(n+1)j\,u^r x^n b^{j-1},
&&j\ge1,\\
L(u^r)&=0.
\end{aligned}
\eea
The nonzero images are nonzero multiples of distinct basis
monomials. They comprise exactly the monomials with $r\ge1$ and
$0\le i\le n-1$, together with all monomials with $i=n$.
Thus
\bea
\ker L=\Q[u],
\quad
R/L(R)\cong
\bigoplus_{\substack{0\le i\le n-1\\j\ge0}}
\Q\,x^i b^j.
\eea

It follows that
\bea
\begin{aligned}
P^{S^1}(\Lambda_gM;\Q)(t)
=\frac1{1-t^2}
+\sum_{i=0}^{n-1}\sum_{j\ge0}t^{d-1+id+j\ell}
=\frac1{1-t^2}
+\frac{t^{d-1}}{1-t^{d(n+1)-2}}
\frac{1-t^{dn}}{1-t^d},
\end{aligned}
\eea
which is (\ref{3.12}).

For the ordinary cohomology, set $u=0$ in the free-loop model and
apply the same finite-filtration argument directly to the resulting
complex. Since $x^{n+1}$ is a non-zero-divisor in $\Q[x,b]$, this gives
the quasi-isomorphic model
\bea
\mathcal C_0=
\frac{\Q[x,b]}{(x^{n+1})}\otimes\Lambda(a),
\quad D_0x=D_0a=0,\quad D_0b=-(n+1)x^na.
\eea
Its cohomology has basis
\bea
1;\quad
[x^ib^j]\quad(1\le i\le n,\ j\ge0);\quad
[x^iab^j]\quad(0\le i\le n-1,\ j\ge0).
\eea
Indeed, the nonclosed even monomials are precisely $b^j$ with
$j\ge1$, whose differentials are linearly independent, and the odd
boundaries are precisely the span of $x^na b^j$, $j\ge0$.
The multiplication in $\mathcal C_0$ therefore yields
\bea
H^*(\Lambda_gM; \Q)
\cong\Q\cdot1\oplus
\left(
\frac{\Lambda^+(x,a)}{(x^{n+1},x^na)}
\otimes\Q[b]
\right).
\eea
Thus the positive-degree expression in part (iv) gives
$\widetilde H^*(\Lambda_gM)$; the displayed summand $\Q\cdot1$
is required for the full cohomology. The ordinary models in
(ii)--(iv) have bounded Betti numbers. Together with the boundedness
criterion, this completes the classification and the proof. \hfill\hb

\medskip

We have the following version of the Morse inequality.

{\bf Theorem 3.3.} (Theorem 6.1 of \cite{Rad2}) {\it Suppose that there exist
only finitely many prime non-contractible closed geodesics $\{c_j\}_{1\le j\le p}$ on Finsler manifold $(M, F)$ with $\pi_1(M)=\Z_2$, and $0\leq a<b\leq\infty$ are regular values of the energy functional $E$.
Define for each $q\in\Z$,
\bea
{M}_q(b,a)&=&
\sum_{1\le j\le p,\;a<E(c_j^{2m-1})<b,\;m\in\N}\rank{\ol{C}}_q(E, c^{2m-1}_j ), if~ F ~is~ irreversible,\nn
\\
{M}_q(b,a)&=&\sum_{1\le j\le p,\;a<E(c_j^{2m-1})<b,\;m\in\N}\left(\rank{\ol{C}}_q(E, c^{2m-1}_j )+\rank{\ol{C}}_q(E, c^{1-2m}_j )\right), if~ F~ is~ reversible.\nn\\ \bar{\beta}_q(b,a)&=&\rank{H}_q(\bar{\Lambda}_g^b,\bar{\Lambda}_g^a),\nn
\eea
where $\bar{\Lambda}_g^b={\Lambda}_g^b/S^1, ~{\Lambda}_g^b=\{\gamma\in\Lambda_g M\mid E(\gamma)<b\}$. Then there holds
\bea
M_q(b,a)&-& M_{q-1}(b,a)
    + \cdots +(-1)^{q}M_0 (b,a)\nn\\&\ge& \bar\beta_q (b,a)-\bar\beta_{q-1}(b,a)
   + \cdots + (-1)^{q} \bar\beta_0 (b,a),~if~b<\infty, \lb{3.47}\\
{M}_q (b,a)&\ge& \bar\beta_q (b,a).\lb{3.48}
\eea
When \(b=\infty\), (3.47) also holds provided that \(M_\ell(\infty,a)<\infty\) for every \(0\le\ell\le q\). In particular, this finiteness condition is satisfied if \(\hat i(c_j)>0\) for every \(1\le j\le p\).

In the following we denote $M_q={M}_q (\infty,0)$ and $\bar\beta_q=\bar\beta_q(\infty,0)$.}

The following lemma follows directly from Lemmas 2.10 and 2.11.

{\bf Lemma 3.4.} {\it  Let $c$ be a non-contractible prime closed geodesic on
a bumpy $(M,\, F)$.
Then
\bea i(c^{p+2})-i(c^p)\in 2\Z,
     \qquad\forall p\in 2\N-1.  \lb{3.49}\eea
}

{\bf Definition 3.5.} {\it Let $c$ be a prime non-contractible closed geodesic on $(M,F)$,
where $\pi_1(M)=\Z_2$. Its local Euler characteristic is
\[
 \chi(c^m)=\sum_{q\ge0}(-1)^q\dim\overline C_q(E,c^m).
\]
For counting in $\Lambda_gM$ define the class-restricted average}
\[
 \hat\chi_g(c)=\lim_{N\to\infty}\frac1N
      \sum_{\substack{1\le m\le N\\m\text{ odd}}}\chi(c^m).
\]

{\bf Remark 3.6.} For a bumpy metric, Proposition~3.1 and Lemma~3.4 give
$\chi(c^m)=(-1)^{i(c)}$ for every odd $m$. Hence
\[
 \hat\chi_g(c)=\frac{(-1)^{i(c)}}2.
\]

By the same proof of Theorem 1.1 of \cite{LiX}, we have the following mean index equality.

{\bf Theorem 3.7.} {\it Suppose $(M,F)$ be a bumpy irreversible Finsler manifold, $\pi_1(M)=\Z_2$, and there are finitely
many prime non-contractible closed geodesics	
$c_1,\ldots,c_p$, all with $\hat i(c_j)>0$. If the average
equivariant Betti number $\bar B(\Lambda_gM;\Q)$ exists, then}
\bea
 \sum_{j=1}^p\frac{\hat\chi_g(c_j)}{\hat i(c_j)}
 =\sum_{j=1}^p\frac{(-1)^{i(c_j)}}{2\hat i(c_j)}
 =\bar B(\Lambda_gM;\Q).
\eea

\setcounter{equation}{0}
\section{Proof of the main theorem  }
%{Section 4}

{\bf Proof of Theorem 1.1 for (i) in Proposition 3.2.}

In this case, there are infinitely many non-contractible closed geodesics
on $(M,\,  F)$, since otherwise the finiteness of the number of prime non-contractible closed geodesics, together with the mean index estimates, implies that the Morse numbers $M_q$ are uniformly bounded in sufficiently high degrees $q$,
note that $b_q(\Lambda_g M)\leq\bar \beta_q+\bar \beta_{q-1}\leq M_q+M_{q-1}$ by the Gysin sequence of the circle bundle:\[S^1\rightarrow ES^1\times\Lambda_g M\rightarrow ES^1\times_{S^1}\Lambda_g M \]
it's a contradiction to the unboundedness of the Betti numbers, see also \cite{GKV} for the Riemannian metric case. \hfill\hb

In the following, we prove Theorem 1.1 by contradiction, thus we assume the following condition.

{\bf (C) There is exactly one non-contractible prime closed geodesic $c$ on $(M, \,F)$.}

{\bf Lemma 4.1.} {\it Under the condition (C), we have $M_{2i-1}=0$ for $i\in\N$.}

{\bf Proof.} By Definition,
\bea M_q=\sum_{m\ge 1}\rank\overline C_q(E, \, c^{2m-1}).\lb{4.1}\eea
Since $E(c)$ is minimal, we have $i(c)=0$,  hence the lemma follows directly from
Proposition 3.1 and Lemma 3.4.\hfill\hb

{\bf Proof of Theorem 1.1 for (iv) in Proposition 3.2.}

In this case, we have $\bar \beta_{d-1}=1$. Hence $M_{d-1}\ge \bar \beta_{d-1}=1$
by Theorem 3.3. This contradicts to Lemma 4.1 since $d\in2\N$.\hfill\hb

{\bf Proof of Theorem 1.1 for (ii) in Proposition 3.2.}

We firstly consider the case that $c$ is non-orientable.

By Lemma 2.9, there is a symplectic path $\gamma: [0, 1] \rightarrow Sp(2d-2)$ with $\gamma(0)=I$
and $\gamma(1)=-P_c$ together with $i(c^m)=i_{-1}(\gamma^m)$ for  $m\in 2\N-1$
 and $i(c^m)=i(\gamma^m)$ for $m\in2\N$ when $d\in2\N$; $\gamma: [0, 1] \rightarrow Sp(2d)$ with $\gamma(0)=I$
and $\gamma(1)=N_1(-1,\,1)\diamond (-P_c)$ together with $i(c^m)=i_{-1}(\gamma^m)$ for  $m\in 2\N-1$ and $i(c^m)=i(\gamma^m)$ for $m\in2\N$ when $d\in2\N-1$.
Hence we can combine these two cases together, i.e., there is a symplectic $\gamma$ with $\gamma(0)=I$
 together with $i(c^m)=i_{-1}(\gamma^m)$ for  $m\in 2\N-1$
 and $i(c^m)=i(\gamma^m)$ for $m\in2\N$.

 Note that $\hat i(c)>0$, otherwise $\hat i(c)=0$ and then by Lemma 2.10, we have $\{i(c^{2m-1})\mid m\in\N\}$ is bounded, which together
 with Proposition 3.1 implies $M_q=0$ for $q$ large enough, it contradicts (3.48).

By Theorem 3.7 and Proposition 3.2, we have $\frac{1}{2\hat{i}(c)}=\bar{B}(\Lambda_g M;\Q)=\frac{1}{2}$,
hence $\hat i(c)=1$.

Denote by   $\mu=\sum_{\theta\in(0, 2\pi)}S_{\gamma(1)}^-(e^{\sqrt{-1}\theta})$ and $\alpha_j=\frac{\theta_j}{\pi}$ where $e^{\sqrt{-1}\theta_j}\in\sigma(\gamma(1))$
for $1\le j\le \mu$.  Let $h=1+\mu$ and
\bea  v=\left(\frac{1}{M\hat i(\gamma, 1)},
\frac{\alpha_1}{\hat i(\gamma, 1)}, \frac{\alpha_2}{\hat i(\gamma, 1)},
\dots\frac{\alpha_{\mu}}{\hat i(\gamma, 1)}
\right)\in\R^h.
\lb{4.2}\eea

Thus by Theorem 2.8,
\bea m_1=\left(\left[\frac{N}{M\hat i(\gamma, 1)}\right]+\xi_1\right)M=N,\lb{4.3}\eea
where we use the fact that $\xi_1=0$, $M=1$, and
$N$ satisfies
\bea \min\left\{\{m_1\alpha_j\},\,1-\{m_1\alpha_j\}\right\}
<\delta/2,\lb{4.4}\eea
for $1\le j\le \mu$ and
\bea i(\gamma, 2m_1-1)&=&	
2N-i(\gamma,\,1)-2S^+_{\gamma(1)}(1)\nn\\&=&	
2N-i(\gamma,\,1) ,\lb{4.5}\\
i(\gamma, 2m_1+1)&=&2N+i(\gamma,\, 1).\lb{4.6}\eea

On the other hand,  we apply Theorem 2.8
to the symplectic path $\gamma^2$. Let $\alpha^\prime_{j}=\frac{2\theta_j}{\pi}$ where $e^{\sqrt{-1}\theta_j}\in\sigma(\gamma(1))$
for $1\le j\le \mu$.  Let $h=1+\mu$ and
\bea  v^\prime=\left(\frac{1}{M\hat i(\gamma, 2)},
\frac{\alpha^\prime_1}{\hat i(\gamma, 2)}, \frac{\alpha^\prime_2}{\hat i(\gamma, 2)},
\dots\frac{\alpha^\prime_{\mu}}{\hat i(\gamma, 2)}
\right)\in\R^h.
\lb{4.7}\eea
By Theorem 2.8 again,
\bea m^\prime_1=\left(\left[\frac{N^\prime}{M\hat i(\gamma, 2)}\right]+\xi_1\right)M
=\frac{N^\prime}{2},\lb{4.8}\eea
where we use the fact that $\xi_1=0$, $M=1$ and $\hat i(\gamma, 2)=2\hat i(\gamma, 1)=2$.

Note that by (\ref{4.4}), we can choose $N^\prime=2N$, note that $v=(1,
\alpha_1, \alpha_2,
\dots,\alpha_{\mu})$,  $v^\prime=(1/2,
\alpha_1, \alpha_2,
\dots,\alpha_{\mu})$, then $N^\prime v^\prime=(N,
2N\alpha_1, 2N\alpha_2,
\dots, 2N\alpha_{\mu})$ and
\bea \min\left\{\{m^\prime_1\alpha^\prime_{j}\},\,1-\{m^\prime_1\alpha^\prime_{j}\}\right\}
<\delta,\lb{4.9}\eea
for $1\le j\le \mu$ by (\ref{4.3}).  Hence we have
\bea m_1=N=m_1^\prime.\lb{4.10}\eea

Thus we have
\bea i(\gamma, 4m_1-2)&=& i(\gamma^2, 2m_1^\prime-1)=	
2N^\prime-i(\gamma^2)-2S^+_{\gamma^2(1)}(1)=4N-i(\gamma,\,2),\lb{4.11}\\
 i(\gamma, 4m_1+2)&=& i(\gamma^2, 2m_1^\prime+1)=	
2N^\prime+i(\gamma^2)=4N+i(\gamma,\,2).\lb{4.12}\eea
Subtracting  (\ref{4.5}) and (\ref{4.6}) from (\ref{4.11}) and (\ref{4.12}) respectively,
we have
\bea i(c^{2N-1})&=&i_{-1}(\gamma,\, 2N-1)=i(\gamma, 4N-2)-i(\gamma, 2N-1)
\nn\\&=&4N-i(\gamma,\,2)-(2N-i(\gamma,\,1))
\nn\\&=&2N-i_{-1}(\gamma,\, 1)=2N-i(c),\lb{4.13}\\
i(c^{2N+1})&=&i_{-1}(\gamma,\, 2N+1)=i(\gamma, 4N+2)-i(\gamma, 2N+1)
\nn\\&=&4N+i(\gamma,\,2)-(2N+i(\gamma,\,1))
\nn\\&=&2N+i_{-1}(\gamma,\, 1)=2N+i(c).\lb{4.14}
\eea
Note that by Lemma 4.1,  Propositions 3.2(ii) and Theorem 3.3, we have $M_q=\bar \beta_q$ for $q\in\Z$.
From (\ref{4.13}) and (\ref{4.14}) and $i(c)=0$, we have
\bea 1=\bar \beta_{2N}=M_{2N}=\sum_{m\ge 1}\rank\overline C_{2N}(E, \, c^{2m-1})\ge 2.\nn\eea
This proves Theorem 1.1 for the case (ii) of Proposition 3.2 when $c$ is non-orientable.

For the orientable case, by Lemma 2.11 and (\ref{4.5})-(\ref{4.6}), we still obtain
$i(c^{2N-1})=2N-i(c)$ and $i(c^{2N+1})=2N+i(c)$, then we finish the proofs as above.
\hfill\hb

Next we consider Theorem 1.1 for (iii) in Proposition 3.2. Firstly we focus on non-orientable case.

Since the Finsler metric $F$ is bumpy, the decomposition in Lemma 2.10 of the minus linearized Poincar\'{e}
map $-P_c$ has the following form:
\bea f_c(1)&=&R(\th_1)\diamond\cdots\diamond R(\th_{r^\prime})
\diamond R(\theta_{r^\prime+1})\diamond\cdots\diamond R(\th_{r})
\nn\\
 &&\diamond N_2(e^{\aa_{1}\sqrt{-1}}, A_{1})\diamond\cdots\diamond N_2(e^{\aa_{r_{\ast}}\sqrt{-1}},A_{r_{\ast}})\nn\\&&\diamond
      N_2(e^{\bb_{1}\sqrt{-1}},B_{1})\diamond\cdots\diamond\,N_2(e^{\bb_{r_{0}}\sqrt{-1}},B_{r_{0}})\diamond H(\pm2)^{\diamond h}, \nn\eea
with  $\theta_j\in(\pi,2\pi)$ for $1\le j\le r^\prime$, $\theta_j\in(0,\pi)$ for $r^\prime< j\le r$  and
$\frac{\th_{j}}{\pi}\not\in\Q$ for $1\le j\le r$; $\frac{\aa_{j}}{2\pi}\not\in\Q$ for $1\le j\le r_{\ast}$;
$\frac{\bb_{j}}{2\pi}\not\in\Q$ for $1\le j\le r_0$ and
\be r + h + 2r_{\ast} + 2r_0 = \dim M-1\equiv n-1, \lb{4.15}\ee
Due to $i(c)=0$, by Lemma 2.10 we obtain
\bea &&i(c^{2m-1})\nn\\
&=&-2(2m-1)r'+2\sum_{j=1}^r \mathcal{E}\left({(2m-1)\theta_{j}\over2\pi}-{1\over2}\right)\nn\\
&=&-2(2m-1)r'+2\sum_{j=1}^r \left(\left({(2m-1)\theta_{j}\over2\pi}-{1\over2}\right)-\left\{{(2m-1)\theta_{j}\over2\pi}-{1\over2}\right\}+1\right)\nn\\
&=&(2m-1)\left(-2r'+\sum_{j=1}^{r}{\theta_{j}\over\pi}\right)-2\sum_{j=1}^r \left(\left\{{(2m-1)\theta_{j}\over2\pi}-{1\over2}\right\}-{1\over2}\right)\nn\\
&=&(2m-1)\left(-2r'+\sum_{j=1}^{r}{\theta_{j}\over\pi}\right)+2r-2\sum_{j=1}^r \left(\left\{{(2m-1)\theta_{j}\over\pi}\right\}+\left\{{-(2m-1)\theta_{j}\over2\pi}\right\}\right),
\label{4.16}\eea
where in the last equality, we have used that
$\{\alpha-\frac{1}{2}\}-\frac{1}{2}=\{2\alpha\}+\{-\alpha\}-1$ for $\alpha\notin\Q$.
Thus
\bea \hat i(c)=-2r^\prime+\sum_{j=1}^r\frac{\theta_j}{\pi}.\lb{4.17}\eea
Let $\hat{\theta}_j=\frac{\th_j}{\pi}-[\frac{\th_j}{\pi}]+1$ for $1\leq j \leq r$ and $\hat{\theta}_j=-\frac{\th_{j-r}}{2\pi}$ for $r+1\leq j \leq 2r$. Then we have
\bea i(c^{2m-1})=(2m-1)\hat i(c)+2r-2\sum_{j=1}^{2r}\left\{(2m-1)\hat\theta_j\right\}.\lb{4.18}
\eea
As in the proof above for the case (ii) of Proposition 3.2, we have $\hat i(c)>0$, then by Theorem 3.7 we have $\frac{1}{2\hat{i}(c)}=\bar{B}(\Lambda_g M;\Q)=\frac{k+1}{2k}$,
hence
\bea \hat i(c)=\frac{k}{k+1}.\lb{4.19}\eea
Then \bea \sum_{j=1}^{2r}\hat\theta_j&=&\sum_{j=1}^{r}\left(\frac{\theta_j}{\pi}-\left[\frac{\theta_j}{\pi}\right]+1\right)+\sum_{j=1}^{r}\left(\frac{-\theta_j}{2\pi}\right)\nn\\
&=&\sum_{j=1}^{r}\left(\frac{\theta_j}{2\pi}-\left[\frac{\theta_j}{\pi}\right]+1\right)\nn\\
&=&\sum_{j=1}^{r}\frac{\theta_j}{2\pi}-r^\prime+r=r+\frac{\hat i(c)}{2}=r+\frac{k}{2(k+1)}\lb{4.20},\eea
where the last equality follows from (\ref{4.17}) and (\ref{4.19}).

Let $s=2r$, then we get the following

{\bf Lemma 4.2. } {\it For non-orientable closed geodesic $c$, the index formulae holds
\bea i(c^{2m-1})&=&(2m-1)\frac{k}{k+1}+s-2\sum_{j=1}^{s}\left\{(2m-1)\hat\theta_j\right\},\lb{4.21}\\
\sum_{j=1}^{s}\hat\theta_j&=&\frac{1}{2}\left(s+\frac{k}{k+1}\right).\lb{4.22}\eea
}

Now we deduce a similar index formula for orientable closed geodesic $c$.

Due to $i(c)=0$, by Lemma 2.11, we obtain\bea i(c^{2m-1})
&=& (2m-1)(0-r )  -r + 2\sum_{j=1}^r{E}\left(\frac{(2m-1)\th_j}{2\pi}\right)\nn\\&=& -(2m-1)r  +r + \sum_{j=1}^r \frac{(2m-1)\th_j}{\pi}-2\sum_{j=1}^r \left\{\frac{(2m-1)\th_j}{2\pi}\right\},\nn\\
\hat i(c)&=&-r+\sum_{j=1}^r\frac{\theta_j}{\pi}.
\nn\eea
Note that (\ref{4.19}) still holds, then we get the similar formulae as in Lemma 4.2:

{\bf Lemma 4.2'. } {\it For orientable closed geodesic $c$, the index formula holds
\bea i(c^{2m-1})
&=& (2m-1)\frac{k}{k+1}  +s  -2\sum_{j=1}^s \left\{(2m-1)\hat\theta_j\right\},\nn\\
\sum_{j=1}^s\hat\theta_j&=&\frac{1}{2}\left(s+\frac{k}{k+1}\right).
\nn\eea
where we set $\hat\theta_j=\frac{\theta_j}{2\pi}, s=r$.}

\medskip

{\bf Lemma 4.3.} {\it Denote by $\dim M=n$, for  any $l\in\N, m\in 2\N-1$,
orientable (or  non-orientable)  closed geodesic $c$,
we always have
\bea |i(c^{m})-2kl|>2n\quad \mbox{holds whenever}\quad  |m-2(k+1)l|>4(k+1)n. \lb{4.23}\eea
}
{\bf Proof.} We only prove it when $c$ is non-orientable since the orientable case is similar. From (\ref{4.21}), we have
\bea  i(c^{m}) = 2kl+\left(m-2(k+1)l\right)\cdot{k\over k+1}+s-2\sum_{j=1}^{s}\left\{{m\hat{\theta}_{j}}\right\}, \lb{4.24}\eea
which, together with $|m-2(k+1)l|>4(k+1)n$, yields immediately that
\bea  |i(c^{m})-2kl|
&\ge& |m-2(k+1)l|\cdot{k\over k+1}-\left|s-2\sum_{j=1}^{s}\left\{{m\hat{\theta}_{j}}\right\}\right|  \nn\\
&>& 4kn-s \ge 4kn-2n\ge 2n,  \nn\eea
where the fact $s\le 2n$ and $k\in\N$ is used. \hfill\hb

Let  $m=2(k+1)l+2L+1$
with $l\in\N$ and $L\in\Z$. Since (\ref{4.21}) and (\ref{4.22}) hold for orientable or non-orientable $c$, then we obtain
\bea  i(c^{m})&=&(2(k+1)l+2L+1)\left(\frac{k}{k+1}\right)+s-2\sum_{j=1}^{s}\left\{(2(k+1)l+2L+1)\hat\theta_j\right\}\nn\\&=& 2kl+s+(2L+1){k\over k+1 }\nn\\
& & -2\left(\left\{m\left(\frac{s}{2}+\frac{k}{2(k+1)}\right)-\sum_{j=2}^{s}\left\{m\hat{\theta}_{j}\right\}\right\}
    +\sum_{j=2}^{s}\left\{m\hat{\theta_{j}}\right\}\right)  \nn\\
&=& 2kl+s+(2L+1){k\over k+1 }\nn\\
& & -2\left(\left\{\frac{s}{2}+\frac{(2L+1)k}{2(k+1)}-\sum_{j=2}^{s}\left\{m\hat{\theta}_{j}\right\}\right\}
    +\sum_{j=2}^{s}\left\{m\hat{\theta_{j}}\right\}\right)  \nn\\
&=& 2kl+2\left[\frac{s}{2}+{(2L+1)k\over2(k+1)}\right]+2\left\{\frac{s}{2}+{(2L+1)k\over2(k+1)}\right\}  \nn\\
& & -2\left(\left\{\left\{{\frac{s}{2}+{(2L+1)k\over2(k+1)}}\right\}-\sum_{j=2}^{s}\left\{m\hat{\theta}_{j}\right\}\right\}
   +\sum_{j=2}^{s}\left\{m\hat{\theta_{j}}\right\}\right)  \nn\\
&=& 2kl+2\left[Q_{L}\right]+2\left\{Q_{L}\right\}-2\left(\left\{\left\{Q_{L}\right\}-\sum_{j=2}^{s}
  \left\{m\hat{\theta}_{j}\right\}\right\}+\sum_{j=2}^{s}\left\{m\hat{\theta_{j}}\right\}\right),
  \label{4.25}\eea
where  $Q_{L}=\frac{s}{2}+{(2L+1)k\over2(k+1)}. $

Since $\sum_{j=2}^{s}\{m\hat{\theta}_{j}\}\in\Q^{c}$, we obtain by (\ref{4.25}) that for $1\leq i\leq s-2$,
\be  i(c^{m})
= \left\{\begin{array}{ll}
2kl+2\left[Q_{L}\right],& \;\mbox{iff}\; \sum_{j=2}^{s}\left\{m\hat{\theta}_{j}\right\}\in (0,\,\left\{Q_{L}\right\}),\\
2kl+2\left[Q_{L}\right]-2i,& \;\mbox{iff}\; \sum_{j=2}^{s}\left\{m\hat{\theta}_{j}\right\}\in (i-1+\left\{Q_{L}\right\},\, i+\left\{Q_{L}\right\}),\\
2kl+2\left[Q_{L}\right]-2(s-1),& \;\mbox{iff}\; \sum_{j=2}^{s}\left\{m\hat{\theta}_{j}\right\}\in (s-2+\left\{Q_{L}\right\},\, s-1).
\end{array}\right.   \lb{4.26}\ee
Let
\bea &&I_{0}(L)=(0,\,\left\{Q_{L}\right\}),\;I_{s-1}(L)=(s-2+\left\{Q_{L}\right\},\,s-1),\;
\nn\\&&I_{i}(L)=(i-1+\left\{Q_{L}\right\},\, i+\left\{Q_{L}\right\})\quad \;for\;\quad 1\le i\le s-2. \lb{4.27}\eea
Then, (\ref{4.26}) can be stated in short as that for any integers $m=2(k+1)l+2L+1$ and  $0\leq i\leq s-1,$
\be   i(c^{m}) = 2kl+2\left[Q_{L}\right]-2i\quad \iff \quad
    \sum_{j=2}^{s}\left\{m\hat{\theta}_{j}\right\}\in I_{i}(L).   \lb{4.28}\ee

Let $\aa=\{\alpha_{1}, \alpha_{2}, \ldots, \alpha_{m}\}$ be a set of $m$ irrational numbers. As usual,
we have

{\bf Definition 4.4.} {\it
The set $\aa$ of irrational numbers is linearly independent over $\Q$, if there do not exist $c_1$,
$c_2$, $\ldots$, $c_m$ in $\Q$ such that $\sum_{j=1}^{m}|c_{j}|> 0$ and
\be  \sum_{j=1}^{m}c_{j}\alpha_{j}\in\Q,  \label{81a}\ee
and is linearly dependent over $\Q$ otherwise. The rank of $\aa$ is defined to be the number of
elements in a maximal linearly independent subset of $\aa$, which we denote by $\rank(\aa)$.
}

{\bf Theorem 4.5.} (Kronecker's approximation theorem, cf. Theorem 7.10 in \cite{Apo}) {\it If $\theta_{1}$, $\theta_{2}$, \dots, $\theta_{\cal R}$ are linearly independent over $\Q$, then the set $\{\{(m\theta_{1}\},\{m\theta_{2}\},\dots,\{m\theta_{\cal R}\})\mid m\in N\}$ is dense in
$$[0,1]^{\cal R}=\underbrace{[0,1]\times[0,1]\times\cdots\times[0,1]}_{\cal R}.$$}

{\bf Lemma 4.6.}  (cf. Lemma 4.1 of \cite{LiX}) {\it
Let $\bar r=\rank(\aa)$. Then there exist $p_{jl}\in\Z$, $\beta_{l}\in\Q^{c}$ and $\xi_{j}\in \Q$ for $1\leq j\leq m$ and $1\leq l\leq \bar r$ such that
 \begin{equation}
\alpha_{j}=\sum_{l=1}^{\bar r}p_{jl}\beta_{l}+\xi_{j},\ \forall 1\leq j\leq m.\label{4.30}\end{equation}
}

Firstly we consider the case that $\rank\{\hat\theta_1,\,\hat\theta_2,\ldots,\hat\theta_{s}\}=1$.  Then by Lemma 4.6, it can be reduced to the following system
\be   \hat{\theta}_{j}=p_{j}\theta+\xi_{j},\ \forall 1\leq j\leq s,  \label{4.31}\ee
with $\theta\in\mathbb{Q}^{c}$, $p_{j}\in \Z\backslash\{0\}$, $\xi_{j}\in\mathbb{Q}\cap[0,1)$
satisfying
\bea
&& p_{1}+p_{2}+\cdots+p_{s}=0, \lb{4.32}\\
&& \{\xi_{1}+\xi_{2}+\cdots+\xi_{s}\}\in(0,1)\backslash\{{1/2}\}, \lb{4.33}\eea
where to get $\xi_j\in [0,1)$, if necessary, we can replace $\hat{\th}_j$ and $\xi_j$ by
$\breve{\theta}_{j}=\hat{\theta}_{j}-[\xi_{j}]$ and $\breve{\xi}_{j}=\{\xi_{j}\}.$
In fact,  by (\ref{4.22}), we have
\bea\frac{1}{2}\left(s+\frac{k}{k+1}\right)=\sum_{j=1}^{s}\hat\theta_j
=\left(\sum_{j=1}^{s}p_j\right)\theta+ \sum_{j=1}^{s}\xi_j,\lb{4.34}\eea
thus (\ref{4.32}) holds since $\theta\in\Q^{c}$, $p_{j}\in \Z\backslash\{0\}$, $\xi_{j}\in\Q$;
(\ref{4.33}) holds since
\bea\{\xi_{1}+\xi_{2}+\cdots+\xi_{s}\}=\left\{\frac{1}{2}\left(s+\frac{k}{k+1}\right)\right\}.\lb{4.35}\eea
Note that
\bea \breve\theta_j=\hat\theta_j-[\xi_j]=p_j\theta+\xi_j-[\xi_j]=p_j\theta+\breve\xi_j.
\eea
Take arbitrarily $\eta\in\Q$ and make  the following natural $\eta$-action  to the system (\ref{4.31}):
\begin{equation}
{\eta}(\theta)=\theta+\eta,\;\ {\eta}(\hat{\theta}_{j})=\hat{\theta}_{j}-\left[\xi_{j}-p_{j}\eta\right],\; {\eta}({\xi}_{j})=\left\{\xi_{j}-p_{j}\eta\right\},\ \forall 1\leq j\leq s,
\lb{4.37}\end{equation}
which is obviously induced by the transformation ${\eta}(\theta)=\theta+\eta$.
In fact,
\bea &&p_j\eta(\theta)+\eta(\xi_j)=p_j(\theta+\eta)+\{\xi_j-p_j\eta\}
=\hat\theta_j-\xi_j+p_j\eta+\{\xi_j-p_j\eta\}\nn\\
=&&\hat\theta_j-((\xi_j-p_j\eta)-\{\xi_j-p_j\eta\})=\hat\theta_j-[\xi_j-p_j\eta]=\eta(\hat\theta_j).
\lb{4.38}\eea
Thus we get a new system
\begin{equation}
{\eta}(\hat{\theta}_{j})=p_{j}{\eta}(\theta)+{\eta}(\xi_{j}),\ \forall 1\leq j\leq s,
\lb{4.39}\end{equation}
with
\bea
\label{4.40}
\{{\eta}({\xi}_{1})+{\eta}({\xi}_{2})+\cdots+{\eta}({\xi}_{s})\}&=&\{\{\xi_{1}-p_{1}\eta\}+\{\xi_{2}-p_{2}\eta\}+\cdots+\{\xi_{s}-p_{s}\eta\}\}\nn\\
&=&\{\xi_{1}+\xi_{2}+\cdots+\xi_{s}-(p_{1}+p_{2}+\dots+p_{s})\eta\}\nn\\
&=&\{\xi_{1}+\xi_{2}+\cdots+\xi_{s}\},
\lb{4.40}\eea
where in the first equality, we have used
$\{\{a\}+\{b\}\}=\{a+b\}$ for $a, b\in\R$; in the third equality we have used the condition (\ref{4.32}).
For simplicity of writing, we also denote  the new system (\ref{4.39})  by (\ref{4.31})$_{\eta}$ meaning that it comes from
(\ref{4.31}) by an $\eta$-action.

For the system (\ref{4.31})$_{\eta}$ with $\eta\in\Q$, we divide the set $\{1\leq j\leq s\}$ into the following three parts:
\bea \mathcal{K}_{0}^{+}(\eta)&=&\{1\leq j\leq s\mid \eta(\xi_{j})=0,\ p_{j}>0\}, \nn\\
\mathcal{K}_{0}^{-}(\eta)&=&\{1\leq j\leq s\mid \eta(\xi_{j})=0,\ p_{j}<0\}, \nn\\
\mathcal{K}_{1}(\eta)&=&\{1\leq j\leq s\mid \eta(\xi_{j})\neq0 \}.\lb{4.41}\eea
Denote by $k_{0}^{+}(\eta)$, $k_{0}^{-}(\eta)$ and $k_{1}(\eta)$ the numbers ${}^{\#}\mathcal{K}_{0}^{+}(\eta)$, ${}^{\#}\mathcal{K}_{0}^{-}(\eta)$ and ${}^{\#}\mathcal{K}_{1}(\eta)$ respectively. For the case of $\eta=0$, we write them for short as $k_{0}^{+}$, $k_{0}^{-}$ and $k_{1}$. It follows immediately that
\bea k_{0}^{+}(\eta)+k_{0}^{-}(\eta)+k_{1}(\eta)=s.\lb{4.42}\eea
By (\ref{4.33}) and (\ref{4.40}), it is obvious that $k_{1}(\eta)\geq1$ for every $\eta\in\Q$.

{\bf Definition 4.7.} {\it For every $\eta\in\Q$, the {\bf absolute difference number} of (\ref{4.31})$_{\eta}$ or (\ref{4.39}) is defined to be the non-negative number $|k_{0}^{+}(\eta)-k_{0}^{-}(\eta)|.$
 The {\bf effective difference number} of {(\ref{4.31})} is defined by
\bea\max\{|k_{0}^{+}(\eta)-k_{0}^{-}(\eta)|\mid \eta\in\Q\}.\lb{4.43}\eea
Two  systems of irrational numbers with rank $1$ are called to be {\bf equivalent}, if their effective difference numbers are the same.
}

{\bf Remark 4.8.} { By the definition of an $\eta$-action in (\ref{4.37}), it can be checked directly that
 $\eta_{1}\circ\eta_{2}=\eta_{1}+\eta_{2}$  for every $\eta_{1}$ and $\eta_{2}$ in $\Q.$ So every system of irrational numbers with rank $1$ is equivalent to the one which comes from itself  by an $\eta$-action.
}

{\bf Theorem 4.9.}
 (cf. Theorem 4.1 of \cite{LiX}) {\it
For every system of irrational numbers (\ref{4.31}) satisfying the conditions (\ref{4.32}) and (\ref{4.33}),  it holds that
   \begin{equation}
\max\{|k_{0}^{+}(\eta)-k_{0}^{-}(\eta)|\mid \eta\in\Q\}\geq1.\lb{4.44}
   \end{equation}
}

By Theorem 4.9 and a permutation of $\{1, 2,\ldots, s\}$ if necessary, we can assume without loss of generality that
\bea |k_{0}^{+}-k_{0}^{-}|\geq1\quad \mbox{and}\quad \mathcal{K}_{1}=\{1,2,\dots, k_{1}\},\lb{4.45}\eea
with $k_{1}\geq1$ due to (\ref{4.33}), and denote by $\xi_{j}={r_{j}\over q_{j}}$ for $1\leq j\leq k_{1}$.

Let $\bar{q}=q_{1}q_{2}\cdots q_{k_1}$ and $m_{l}=2(k+1)\bar{q}l+1$ with $l\in\N$. Then by (\ref{4.31}) we have
\bea
\sum_{j=2}^{s}\left\{m_{l}\hat{\theta}_{j}\right\}&=&\sum_{j=2}^{k_1}\left\{m_{l}\hat{\theta}_{j}\right\}+\sum_{j=k_{1}+1}^{s}\left\{m_{l}\hat{\theta}_{j}\right\}\nn\\
&=&\sum_{j=2}^{k_1}\left\{m_{l}(p_j\theta+\xi_j)\right\}+\sum_{j=k_{1}+1}^{s}\left\{m_{l}p_j\theta\right\}\nn\\&=&\sum_{j=2}^{k_1}\left\{p_{j}\{m_{l}\theta\}+\xi_{j}\right\}+\sum_{j=k_{1}+1}^{s}\left\{p_{j}\{m_{l}\theta\}\right\},\lb{4.46}
\eea
where in the second equality, we have used $\xi_j=0$ for $k_1+1\le j\le s$;
in the third equality, we have used $q_j|(m_l-1)$ and $\{\{a\}+\{b\}\}=\{a+b\}$.

Since $\theta\in\mathbb{Q}^{c}$,  the set $\{\{m_{l}\theta\}\mid l\in\N\}$ is dense in $[0,1]$
by Theorem 4.5. For every $L\in\Z$, define the  auxiliary function
\begin{equation}
f_{L}(x)=\sum_{j=2}^{k_1}\left\{\left\{p_{j}x+\xi_{j}\right\}+2L\hat{\theta}_{j}\right\}+\sum_{j=k_{1}+1}^{s}\left\{\left\{p_{j}x\right\}+2L\hat{\theta}_{j}\right\},\ \forall x\in [0,1],
\lb{4.47}\end{equation}
and denote for simplicity by $f=f_{0}$.
Note that the discontinuous points of $f_L$ are those $x$ satisfying
$p_jx+\xi_j\in\Z$, $\left\{p_{j}x+\xi_{j}\right\}+2L\hat{\theta}_{j}\in\Z$,  $p_jx\in\Z$
and $\left\{p_{j}x\right\}+2L\hat{\theta}_{j}\in\Z$, thus there are only finitely many.

Let $a$ and $b$ in $(0,1)$ be two real numbers  sufficiently close to $0$ and $1$ respectively.
Then,
\bea
f(a)&=&\sum_{j=2}^{k_1}\left\{p_{j}a+\xi_{j}\right\}+\sum_{j=k_{1}+1}^{s}\left\{p_{j}a\right\}
\nn\\&=&\sum_{j=2}^{k_1}\left(p_{j}a+\xi_{j}\right)+\sum_{j\in\mathcal{K}_{0}^{+}}p_{j}a+\sum_{j\in\mathcal{K}_{0}^{-}}\left(1+p_{j}a\right)\nn\\
&=&k_{0}^{-}+\sum_{j=2}^{s}p_{j}a+\sum_{j=2}^{k_1}\xi_{j},
\lb{4.48}\eea
where in the second equality, we have used
$\xi_j\in(0,\,1)\cap\Q$ and $a\rightarrow 0+$,
hence $\{p_ja+\xi_{j}\}=p_ja+\xi_j$ for $2\le j\le k_1$;
$p_j>0$ for $j\in\mathcal{K}_0^+$, hence $\{p_ja\}=p_ja$ for $j\in\mathcal{K}_0^+$;
$p_j<0$ for $j\in\mathcal{K}_0^-$, hence $\{p_ja\}=1+p_ja$ for $j\in\mathcal{K}_0^-$.
\bea
f(b)&=&\sum_{j=2}^{k_1}\left\{p_{j}b+\xi_{j}\right\}+\sum_{j=k_{1}+1}^{s}\left\{p_{j}b\right\}
\nn\\&=&\sum_{j=2}^{k_1}\left\{p_{j}(b-1)+\xi_{j}\right\}+\sum_{j=k_{1}+1}^{s}\left\{p_{j}(b-1)\right\}
\nn\\&=&\sum_{j=2}^{k_1}\left(p_{j}(b-1)+\xi_{j}\right)+\sum_{j\in\mathcal{K}_{0}^{+}}(p_{j}(b-1)+1)+\sum_{j\in\mathcal{K}_{0}^{-}}p_{j}(b-1)\nn\\
&=&k_{0}^{+}+\sum_{j=2}^{s}p_{j}(b-1)+\sum_{j=2}^{k_1}\xi_{j},
\lb{4.49}\eea
where in the third equality, we have used
$\xi_j\in(0,\,1)\cap\Q$ and $b\rightarrow 1-$,
hence $\{p_j(b-1)+\xi_{j}\}=p_j(b-1)+\xi_j$;
 for $2\le j\le k_1$;
$p_j>0$ for $j\in\mathcal{K}_0^+$, hence $\{p_j(b-1)\}=p_j(b-1)+1$ for $j\in\mathcal{K}_0^+$;
$p_j<0$ for $j\in\mathcal{K}_0^-$, hence $\{p_j(b-1)\}=p_j(b-1)$ for $j\in\mathcal{K}_0^-$.

It follows by (\ref{4.48}) and (\ref{4.49}) that
\begin{equation}
\label{4.50}
|f(b)-f(a)|=\left|k_{0}^{+}-k_{0}^{-}+\sum_{j=2}^{s}p_{j}(b-1-a)\right|=|k_{0}^{+}-k_{0}^{-}+p_{1}(-b+1+a)|,
\end{equation}
where the second identity we have used $\sum_{j=1}^{s}p_{j}=0$.

{\bf Lemma 4.10.}  {\it Given $\bar{N}\in\N$, for any $a$ and $b$ in $(0,1)$ sufficiently close to $0$ and $1$ respectively, we have

{\rm (i)} $f(a)$ and $f(b)$ lie in different intervals of  (\ref{4.27}) with $L=0$,

{\rm (ii)} For any $1\leq |L|\leq\bar{N}$, $f_{L}(a)$,  $f_{L}(b)$ and $f_{L}(0)$ lie in the same interval of (\ref{4.27}).
}

{\bf Proof.}  (i) By (\ref{4.50}) and $a\rightarrow 0+$, $b\rightarrow 1-$, $|f(b)-f(a)|\approx|k_{0}^{+}-k_{0}^{-}|.$
 Here and later,  we use $X\approx Y$ as usual to mean that $X$ is sufficiently close to $Y$ in the context of writing.
 Since the length of each interval in (\ref{4.27}) with $L=0$ is less than or equal to $1$, so $f(a)$ and $f(b)$
 must lie in different ones, provided that  $|k_{0}^{+}-k_{0}^{-}|\geq2$.

 If $|k_{0}^{+}-k_{0}^{-}|=1$,  then $|f(b)-f(a)|\approx1$.
 Note that we have $s\ge 2$. We have the following three subcases.

{\bf Subcase I:}\quad  If $k_{1}\geq2$,
\bea\left\{k_{0}^{-}+\sum_{j=2}^{k_1}\xi_{j}\right\}&=&\left\{-\xi_{1}+\sum_{j=1}^{k_1}\xi_{j}\right\}=\left\{\sum_{j=1}^{s}\xi_{j}-\xi_{1}\right\}\nn\\&=&\left\{\left\{\sum_{j=1}^{s}\xi_{j}\right\}-\xi_{1}\right\}=\left\{ \{Q_{0}\}-\xi_{1}\right\},\lb{4.51}\eea
where in the second equality, we have used $\xi_j=0$ for $k_1+1\le j\le s$;
in the last equality, we have used (\ref{4.35}) and (\ref{4.25}), i.e., $\{Q_0\}=\left\{\frac{s}{2}+\frac{k}{2(k+1)}\right\}$.

Notice that the dividing points of the intervals in (\ref{4.27}) with $L=0$ are $$0,\ \{Q_{0}\},\ 1+\{Q_{0}\},\ 2+\{Q_{0}\},\ \dots,\ s-2+\{Q_{0}\},\ s-1.$$
Since $\xi_j\in(0, \,1)$ for $1\le j\le k_1$ and $k_1\ge 2$, the point
$k_{0}^{-}+\sum_{j=2}^{k_1}\xi_{j}$ must be an interior point of these intervals.
In fact, $k_{0}^{-}+\sum_{j=2}^{k_1}\xi_{j}\ge \sum_{j=2}^{k_1}\xi_j>0$;
$k_0^-+\sum_{j=2}^{k_1}\xi_j<k_0^-+k_1-1\le s-1$.
If $k_{0}^{-}+\sum_{j=2}^{k_1}\xi_{j}=k+\{Q_0\}$ for some $k\in\{0, 1,\ldots, s-2\}$,
then $\{k_{0}^{-}+\sum_{j=2}^{k_1}\xi_{j}\}=\{Q_0\}$, which contradicts  (\ref{4.51}),
i.e.,  $\left\{ \{Q_{0}\}-\xi_{1}\right\}\neq\{Q_0\}$ since $\xi_1\in(0,\,1)$.
Hence by (\ref{4.48}),  $f(a)$ is an interior point of the intervals in (\ref{4.27}) with $L=0$.
Note that $|f(b)-f(a)|\approx1$, this  yields that $f(a)$ and $f(b)$ must lie in two different intervals
since each interval has length less or equal to $1$.

{\bf Subcase II:} If $k_{1}=1$ and $k_{0}^{-}\geq1$, then by (\ref{4.32}) and (\ref{4.48}), we get $f(a)=k_{0}^{-}-p_{1}a\neq k_{0}^{-}, f(a)\approx k_{0}^{-}$ is also an interior point and (i) follows.

{\bf Subcase III:} If $k_{1}=1$ and $k_{0}^{-}=0$, then $f(a)=\sum_{j=2}^{s}p_{j}a$ lies in the first interval whose length is $\{Q_{0}\}<1$ and so $f(b)$ must lie in another one.

(ii) It can be checked directly that $\lim_{a\to0}f_{L}(a)=\lim_{b\to1}f_{L}(1)=f_{L}(0)\in\mathbb{Q}^{c}$ since $\xi_{j}\in\Q$ for $1\leq j\leq s$ and $\sum_{j=2}^{s}2L\hat{\theta}_{j}\in\mathbb{Q}^{c}$.
In fact, for $\mu\in (0,\,1)$ fixed, clearly $\{\mu\pm \epsilon\}=\mu\pm\epsilon\rightarrow\mu=\{\mu\}$ for $\epsilon\rightarrow 0$.
This together with $\xi_j\in\Q$ and $2L\hat\theta_j\in\Q^c$  implies  for $x\rightarrow 0+, 1-$,
\bea f_{L}(x)&=&\sum_{j=2}^{k_1}\left\{\left\{p_{j}x+\xi_{j}\right\}+2L\hat{\theta}_{j}\right\}+\sum_{j=k_{1}+1}^{s}\left\{\left\{p_{j}x\right\}+2L\hat{\theta}_{j}\right\}
\nn\\&\rightarrow&\sum_{j=2}^{k_1}\left\{\left\{\xi_{j}\right\}+2L\hat{\theta}_{j}\right\}+\sum_{j=k_{1}+1}^{s}\left\{2L\hat{\theta}_{j}\right\}=f_L(0).\lb{4.52}\eea
Note that
\bea \{f_L(0)\}&=&\left\{\sum_{j=2}^{k_1}\left\{\left\{\xi_{j}\right\}+2L\hat{\theta}_{j}\right\}+\sum_{j=k_{1}+1}^{s}\left\{2L\hat{\theta}_{j}\right\}\right\}
\nn\\&=&\left\{\sum_{j=2}^{s}(\xi_{j}+2L\hat{\theta}_{j})\right\}
=\left\{-2Lp_1\theta+(2L+1)\sum_{j=2}^{s}\xi_{j}\right\}\in\Q^c,
\eea
where the last equality follows from (\ref{4.32}).

But the dividing points of these intervals in (\ref{4.27}) with $1\leq|L|\leq\bar{N}$ are finitely rational numbers, so $f_{L}(0)$ is an interior point of these intervals and (ii) follows. \hfill\hb

\medskip

{\bf Proof of Theorem 1.1 for (iii) in Proposition 3.2 with $\rank\{\hat\theta_1,\ldots,\hat\theta_{s}\}=1$.}

Note that $f$ contains only finitely many discontinuous points $(0,1)$. Without loss of generality, we assume $a$ and $b$ to be two continuous points of  $f$ and choose $l_{1}$, $l_{2}\in\N$ with $l_{2}-l_{1}$ sufficiently large such that $\{m_{l_1}\theta\}\approx a$ and $\{m_{l_2}\theta\}\approx b$, this follows directly from Theorem 4.5. Then by (\ref{4.46}), (\ref{4.47}) and (i) of Lemma 4.10, we get
$\sum_{j=2}^{s}\left\{m_{l_1}\hat{\theta}_{j}\right\}$ and $\sum_{j=2}^{s}\left\{m_{l_2}\hat{\theta}_{j}\right\}$ lie in different intervals of (\ref{4.27}) with $L=0$.  Suppose that
\bea\sum_{j=2}^{s}\left\{m_{l_1}\hat{\theta}_{j}\right\}\in I_{i^\prime}\;and\; \sum_{j=2}^{s}\left\{m_{l_2}\hat{\theta}_{j}\right\}\in I_{i^{\prime\prime}},\lb{4.54}\eea
with $\{i^\prime,\, i^{\prime\prime}\}\subseteq\{0,1,2,\dots,s-1\}$ and $i^{\prime}\neq i^{\prime\prime}$. By (\ref{4.28}) we have
\bea i(c^{m_{l_{1}}})&=&2k\bar{q}l_{1}+2[Q_{0}]-2i^{\prime},\nn\\
i(c^{m_{l_2}})&=&2k\bar{q}l_{2}+2[Q_{0}]-2i^{\prime\prime}.\lb{4.55}
\eea

Note that by (\ref{3.10}), we have
\bea \bar \beta_{2k\bar{q}l_{1}+2[Q_{0}]-2i^{\prime\prime}}=\bar \beta_{2k\bar{q}l_{2}+2[Q_{0}]-2i^{\prime\prime}}\equiv\beta.\lb{4.56}\eea
By (\ref{4.27}),
\bea \left|2[Q_{0}]-2i^{\prime\prime}\right|&=&\left|2\left[\frac{s}{2}+\frac{k}{2(k+1)}\right]-2i^{\prime\prime}\right|\nn\\
&\le& \max\left\{2\left[\frac{s}{2}+\frac{k}{2(k+1)}\right],\;
2(s-1)-2\left[\frac{s}{2}+\frac{k}{2(k+1)}\right]\right\}\nn\\
&\le&s\le 2n.\lb{4.57}\eea
Take $\bar N>4(k+1)n$ in (ii) of Lemma 4.10, we have the following

{\bf Claim.} {\it  There exist $L_{i}\in\Z$ with $1\leq |L_{i}|\leq \bar{N}$ and  $1\leq i\leq\beta$ such that
\bea i(c^{m_{l_{1}}+2L_{i}})=2k\bar{q}l_{1}+2[Q_{0}]-2i^{\prime\prime}.\lb{4.58}\eea
}

In fact, by (\ref{4.56}) and Theorem 3.3, there exist $\bar m_i$ for $1\le i\le \beta$
such that $i(c^{\bar m_i})=2k\bar{q}l_{1}+2[Q_{0}]-2i^{\prime\prime}$.
Write $\bar m_i=m_{l_1}+2L_i$ with $L_i\in\Z$ for $1\le i\le \beta$.
Then by (\ref{4.57}) and Lemma 4.3, we have $|L_i|\le \bar N$.
Since $i^{\prime}\neq i^{\prime\prime}$, we have $L_i\neq 0$ by (\ref{4.55}).
This proves the claim.

Since $\{m_{l_1}\theta\}\approx a$ and $\{m_{l_2}\theta\}\approx b$,
by (\ref{4.47}) that
\bea \sum_{j=2}^{s}\left\{(m_{l_1}+2L_i)\hat\theta_j\right\}\approx f_{L_i}(a),
\quad \sum_{j=2}^{s}\left\{(m_{l_2}+2L_i)\hat\theta_j\right\}\approx f_{L_i}(b).\lb{4/59}\eea
Thus by (ii) of Lemma 4.10, $\sum_{j=2}^{s}\left\{(m_{l_1}+2L_i)\hat\theta_j\right\}\approx f_{L_i}(a)$ and $\sum_{j=2}^{s}\left\{(m_{l_2}+2L_i)\hat\theta_j\right\}\approx f_{L_i}(b)$.
belong to the same interval of (\ref{4.27}) for $1\le |L_i|\le \bar N$.
Hence by (\ref{4.28}), we have
\begin{equation}
\label{4.60}
i(c^{m_{l_{2}}+2L_{i}})=2k\bar{q}l_{2}+2[Q_{0}]-2i^{\prime\prime},\ \forall 1\leq i\leq \beta.
\end{equation}
By (\ref{4.55}) and (\ref{4.60}), it yields $M_{2k\bar{q}l_{2}+2[Q_{0}]-2i^{\prime\prime}}\ge\beta+1$ which contradicts (\ref{4.56}), since $M_q=\bar \beta_q$ for any $q\in\Z$. \hfill\hb

\medskip

Next we consider the case that  $\mathcal{R}=\rank\{\hat\theta_1,\ldots,\hat\theta_{s}\}\ge 2$.

By Lemma 4.6,  there are $p_{jl}\in\Z$, $\theta_{k_l}\in\Q^{c}$ and $\xi_{j}\in\Q$ with $1\leq l\leq \cal{R}$ and $1\leq j\leq s$ such that
\begin{equation}
\label{4.61}
\hat{\theta}_{j}=\sum_{l=1}^{\cal R}p_{jl}\theta_{k_l}+\xi_{j},\ \forall 1\leq j\leq s.
\end{equation}
Moreover,  $\theta_{k_1}$, $\theta_{k_2}$ \dots, $\theta_{k_{\cal R}}$ are linearly independent over $\Q$. Due to
(\ref{4.22}), it follows
\begin{equation}
\label{4.62}
\sum_{j=1}^{s}p_{jl}=0,\ \forall 1\leq l\leq \cal R.
\end{equation}
In fact,  by (\ref{4.22}), we have
\bea\frac{1}{2}\left(s+\frac{k}{k+1}\right)=\sum_{j=1}^{s}\hat\theta_j
=\sum_{l=1}^{\cal R}\left(\sum_{j=1}^{s}p_{jl}\right)\theta_{k_l}+ \sum_{j=1}^{s}\xi_j,\lb{4.63}\eea
thus (\ref{4.62}) holds since $\theta_{k_1}$, $\theta_{k_2}$ \dots, $\theta_{k_{\cal R}}$ are linearly independent over $\Q$, $p_{jl}\in \Z$, $\xi_{j}\in\Q$.
Hence
\bea\{\xi_{1}+\xi_{2}+\cdots+\xi_{s}\}=\left\{\frac{1}{2}\left(s+\frac{k}{k+1}\right)\right\}\in (0,\,1)\setminus\{1/2\}.\lb{4.64}\eea

Our basic idea for proving this case  is to construct an irrational system  with rank $1$  associated to (\ref{4.61}),  then use the similar argument as above to prove Theorem 1.1.

{\bf Lemma 4.11.} (cf. Lemma 5.2 of \cite{LiX}) {\it
For the integers $p_{jl}$s in (\ref{4.61}) and (\ref{4.62}),  there are $s_{2},\ s_{3}, \dots,\ s_{\cal R}\in\Z$ such that
\begin{equation}
\label{4.65}
p_{j1}+\sum_{l=2}^{\cal R}s_{l}p_{jl}\in\Z\backslash\{0\},\ \forall 1\leq j\leq s.
\end{equation}
}

Let $\tilde p_{j1}=p_{j1}+\sum_{l=2}^{\cal R}s_lp_{jl}\in\Z\setminus\{0\}$ and
$\tilde p_{jl}=p_{jl}$ for $2\le l\le\cal R$. By Lemma 4.11, we can make the change of variables
$\td{\theta}_{k_1}=\theta_{k_1}$  and $\td{\theta}_{k_l}=\theta_{k_l}-s_{l}\theta_{k_1}$ for $2\leq l\leq \cal R$.
Then the system (\ref{4.61}) is transformed to
\begin{equation}
\label{4.66}
\hat{\theta}_{j}=\sum_{l=1}^{\cal R}\td{p}_{jl}\td{\theta}_{k_l}+\xi_{j},\ \forall 1\leq j\leq s,
\end{equation}
with
$\td{p}_{j1}=p_{j1}+\sum_{l=2}^{\cal R}s_{l}p_{jl}\in\Z\backslash\{0\},$ and by (\ref{4.62}) we have
\bea
\label{4.67}
\sum_{j=1}^{s}\td{p}_{j1}=\sum_{j=1}^{s}p_{j1}+\sum_{j=1}^{s}\sum_{l=2}^{\cal R}s_{l}p_{jl}=0+\sum_{l=2}^{\cal R}s_{l}\left(\sum_{j=1}^{s}p_{jl}\right)=0.
\eea
Since  $\theta_{k_1}$, $\theta_{k_2}$, \dots, $\theta_{k_{\cal R}}$ are linearly independent over $\Q$, so do  $\td{\theta}_{k_1}$, $\td{\theta}_{k_2}$, \dots, $\td{\theta}_{k_{\cal R}}.$

For \(1\le j\le s\), let
\[
 {\bf p}_j=(\widetilde p_{j1},\ldots,\widetilde p_{j{\cal R}})\in\Z^{\cal R},
 \qquad
 \sum_{j=1}^{s}{\bf p}_j={\bf 0},
 \qquad
 \widetilde p_{j1}\ne0.
\]
\[
 T:=\sum_{j=1}^{s}\widehat\theta_j\in\Q,
 \qquad
 \left\{\sum_{j=1}^{s}\xi_j\right\}
 =\{T\}=\{s/2\}+\frac{k}{2(k+1)}.
\]

Consider the following irrational system  with rank $1$ associated to (\ref{4.66})
\begin{equation}
\label{4.68}
\hat{\alpha}_{j}=\td{p}_{j1}\td{\theta}_{k_1}+\xi_{j},\ \forall 1\leq j\leq s.
\end{equation}
By Theorem 4.9, (\ref{4.61}) and (\ref{4.64}), without loss of generality we can assume for (\ref{4.68}) that $|\td{k}_{0}^{+}-\td{k}_{0}^{-}|\geq1$ and
$\td{\mathcal{K}}_{1}=\{1,2,\dots, \td{k}_{1}\}$, and denote by $\xi_{j}={r_{j}\over q_{j}}$ for $1\leq j\leq \td{k}_{1}.$

Then we have
\[
 0<\xi_j<1\quad(1\le j\le\td{k}_{1}),
 \qquad
 \xi_j=0\quad(\td{k}_{1}<j\le s),
 \qquad \td{k}_{1}\ge1.
\]
\[
 \td{k}_{0}^{+}=\#\{j:\xi_j=0,\ \widetilde p_{j1}>0\},
 \qquad
 \td{k}_{0}^{-}=\#\{j:\xi_j=0,\ \widetilde p_{j1}<0\}.
\]

Let $\td{q}=q_{1}q_{2}\cdots q_{\td{k}_{1}}$ and $\td{m}_{l}=2(k+1)\td{q}l+1$ for $l\in\N$. Then, we get by (\ref{4.66}) that
\bea
\sum_{j=2}^{s}\left\{\td{m}_{l}\hat{\theta}_{j}\right\}
&=&\sum_{j=2}^{\td{k}_1}\left\{\td{m}_{l}\hat{\theta}_{j}\right\}+\sum_{j=\td{k}_{1}+1}^{s}\left\{\td{m}_{l}\hat{\theta}_{j}\right\}\nn\\
&=&\sum_{j=2}^{\td{k}_1}\left\{\sum_{l=1}^{\cal R}\td{p}_{jl}\{\td{m}_{l}\td{\theta}_{k_l}\}+\xi_{j}\right\}+\sum_{j=\td{k}_{1}+1}^{s}\left\{\sum_{l=1}^{\cal R}\td{p}_{jl}\{\td{m}_{l}\td{\theta}_{k_l}\}\right\},\lb{4.69}
\eea
where $\xi_j=0$ for $\tilde k_1+1\le j\le s$ is used.
By Theorem 4.5,  the set $\{(\{\td{m}_{l}\td{\theta}_{k_1}\},\{\td{m}_{l}\td{\theta}_{k_2}\},\dots,\{\td{m}_{l}\td{\theta}_{k_{\cal R}}\})\mid l\in\N\}$ is dense in $[0,1]^{\cal R}$. For every $L\in\Z$, we can introduce the  auxiliary multi-variable function on $[0,1]^{\cal R}$,
\bea g_{L}(x_{1},x_{2},\dots,x_{\cal R})=\sum_{j=2}^{\td{k}_1}\left\{\sum_{l=1}^{\cal R}\td{p}_{jl}x_{l}+\xi_{j}+2L\hat{\theta}_{j}\right\}
+\sum_{j=\td{k}_{1}+1}^{s}\left\{\sum_{l=1}^{\cal R}\td{p}_{jl}x_{l}+2L\hat{\theta}_{j}\right\},
\lb{4.70}\eea
and denote for simplicity by $g=g_{0}$.

{\bf Lemma 4.12.}
Assume the notation above, for every \(\overline N\in\N\), there exist constants
\(\eta,\delta\in(0,1)\) such that, whenever
\({\bf a},{\bf b}\in(0,1)^{\cal R}\) satisfy
\[
 \begin{aligned}
 0<a_1<\delta,\qquad&
 0<a_\ell<\eta a_1 &&(2\le\ell\le {\cal R}),\\
 0<1-b_1<\delta,\qquad&
 0<b_\ell<\eta(1-b_1) &&(2\le\ell\le {\cal R}),
 \end{aligned}
\]
the following assertions hold:
\begin{enumerate}
 \renewcommand{\labelenumi}{(\roman{enumi})}
 \item \(g({\bf a})\) and \(g({\bf b})\) belong to different intervals
 in (4.27) with \(L=0\).
 \item For every integer \(L\) with \(1\le|L|\le\overline N\),
 the three numbers \(g_L({\bf a})\), \(g_L({\bf b})\), and \(g_L({\bf 0})\)
 belong to the same interval in (4.27).
\end{enumerate}

{\bf Proof.} Write \({\bf e}_1=(1,0,\ldots,0)\). Choose \(\eta\in(0,1)\)
sufficiently small that
\[
 \eta\sum_{\ell=2}^{\cal R}|\widetilde p_{j\ell}|
 <\frac12|\widetilde p_{j1}|,
 \qquad 1\le j\le s.
\]
For points satisfying the inequalities in the statement, we have
\[
 \left|\sum_{\ell=2}^{\cal R}\widetilde p_{j\ell}a_\ell\right|
 <\frac12|\widetilde p_{j1}|a_1,
 \qquad
 \left|\sum_{\ell=2}^{\cal R}\widetilde p_{j\ell}b_\ell\right|
 <\frac12|\widetilde p_{j1}|(1-b_1).
\]
Thus \({\bf p}_j\cdot{\bf a}\) has the sign of \(\widetilde p_{j1}\),
whereas \({\bf p}_j\cdot({\bf b}-{\bf e}_1)\) has the opposite sign.
Both quantities tend uniformly to zero as \(\delta\to0\).

By taking \(\delta\) sufficiently small, we may therefore ensure
that their absolute values are less than one and that
\[
 \begin{aligned}
 0&<\xi_j+{\bf p}_j\cdot{\bf a}<1,\\
 0&<\xi_j+{\bf p}_j\cdot({\bf b}-{\bf e}_1)<1
 \end{aligned}
 \qquad(1\le j\le\td{k}_{1}).
\]
Set
\[
 S=\sum_{j=2}^{\td{k}_{1}}\xi_j,\qquad
 A=\td{k}_{0}^{-}+S,\qquad B=\td{k}_{0}^{+}+S,
\]
where an empty sum is understood to be zero.
The preceding sign information, together with the integrality of
\(\widetilde p_{j1}\), gives
\[
 \begin{aligned}
 g({\bf a})
 &=\td{k}_{0}^{-}+S+\sum_{j=2}^{s}{\bf p}_j\cdot{\bf a}
   =A-{\bf p}_1\cdot{\bf a},\\
 g({\bf b})
 &=\td{k}_{0}^{+}+S+\sum_{j=2}^{s}{\bf p}_j\cdot({\bf b}-{\bf e}_1)
   =B-{\bf p}_1\cdot({\bf b}-{\bf e}_1).
 \end{aligned}
\]
Consequently, these values converge uniformly to \(A\) and \(B\),
respectively, as \(\delta\to0\).

Let \(q=\{Q_0\}=\{s/2+k/[2(k+1)]\}\in(0,1)\). By the assumptions,
\[
 \{A\}=\{B\}
 =\left\{\sum_{j=1}^{s}\xi_j-\xi_1\right\}
 =\{q-\xi_1\}\ne q,
\]
since \(0<\xi_1<1\). Hence neither \(A\) nor \(B\) is an
internal dividing point of the intervals in (4.27)
with \(L=0\). Moreover,
\[
 A,B\in[0,s-1],\qquad |A-B|=|\td{k}_{0}^{+}-\td{k}_{0}^{-}|\ge1.
\]
Each of \(A\) and \(B\) therefore belongs to the closure of a
unique interval in that partition. These two intervals are
distinct: otherwise, their length, which is at most one, would
force \(A\) and \(B\) to be the two endpoints of an interval of
length one, contradicting their avoidance of the internal
dividing points.

All summands defining \(g({\bf a})\) and \(g({\bf b})\) lie in \((0,1)\),
so
\[
 0<g({\bf a}),g({\bf b})<s-1.
\]
Their uniform convergence to \(A\) and \(B\) now implies that,
for sufficiently small \(\delta\), they lie in the respective
distinct intervals determined above. This also covers the cases
where a limiting value is \(0\) or \(s-1\).
Assertion~(i) follows.

To prove~(ii), fix a nonzero integer \(L\). Since each
\(\widehat\theta_j\) is irrational and each \(\xi_j\) is rational,
\[
 \xi_j+2L\widehat\theta_j\notin\Z.
\]
Thus \(g_L\) is continuous at both \({\bf 0}\) and \({\bf e}_1\).
Furthermore,
\[
 g_L({\bf e}_1)=g_L({\bf 0}),
\]
because every \(\widetilde p_{j1}\) is an integer.

We claim that this common value lies strictly inside an interval
in (4.27). Indeed,
\[
 g_L({\bf 0})=\sum_{j=2}^{s}\{\xi_j+2L\widehat\theta_j\},
\]
and, modulo an integer, this equals
\[
 S+2L\left(\sum_{j=1}^{s}\widehat\theta_j-\widehat\theta_1\right)
 =S+2L(T-\widehat\theta_1).
\]
Since \(S,T\in\Q\) and \(\widehat\theta_1\notin\Q\),
the number \(g_L({\bf 0})\) is irrational.
It therefore avoids every dividing point of
(4.27), all of which are rational.
Also, every summand in its definition lies in \((0,1)\),
so \(g_L({\bf 0})\in(0,s-1)\).

By continuity, there exist neighborhoods of \({\bf 0}\) and \({\bf e}_1\)
whose images under \(g_L\) lie in the same interval containing
\(g_L({\bf 0})\). Since only finitely many integers satisfy
\(1\le|L|\le\overline N\), decreasing \(\delta\) once more makes
this conclusion valid simultaneously for all such \(L\).
This proves~(ii).\hfill\hb

{\bf Proof of Theorem 1.1 for (iii) in Proposition 3.2 with $\mathcal{R}=\rank\{\hat\theta_1,\ldots,\hat\theta_{s}\}\ge 2$.}
We apply Lemma 4.12.
Both conical regions are nonempty open sets, and every tail of the Kronecker sequence is dense. Hence the second iteration parameter can be chosen arbitrarily large while ensuring that the conclusions of Lemma 4.12 hold simultaneously for all shifts in the prescribed finite range. The same counting argument as in the rank-one case then yields a contradiction.
\hfill\hb

\medskip

{\bf Proof of Theorem 1.2.}

Proofs for the cases (i) and (iv) in Proposition 3.2 are the same as those of Theorem 1.1; thus we only need to consider (ii) and (iii) in Proposition 3.2 in the following.

We assume that there exists only one non-contractible closed geodesic $c$ on $(M, \,F)$ with bumpy reversible Finsler metric.  Then its inverse curve $c^{-1}$ has the same energy and plays the same role in the variational setting of the energy functional $E$ on $\Lambda M$ as $c$. In particular, the $m$-th iterates $c^m$ and $c^{-m}$ of a closed geodesic $c$ and its inverse curve $c^{-1}$ have precisely the same Morse indices, nullities and critical modules, but $S^1\cdot c\neq S^1 \cdot c^{-1}$. By Definition,
\bea M_q=\sum_{m\ge 1}\left(\rank\overline C_q(E, \, c^{2m-1})+\rank\overline C_q(E, \, c^{1-2m})\right).\nn\eea
Note that $\hat i(c)>0$, otherwise $\hat i(c)=0$ and then by Lemma 2.10 and Lemma 2.11, we have $\{i(c^{2m-1})\mid m\in\N\}$ is bounded, which together
 with Proposition 3.1 implies $M_q=0$ for $q$ large enough, it contradicts (3.48). Then by the definition of $M_q$ and Proposition 3.1, $M_q$ is finite.

Since $E(c)=E(c^{-1})$ is minimal, we have $i(c)=i(c^{-1})=0$,  hence $M_q=0$  for odd $q$  and $M_0\geq2$ by Proposition 3.1 and Lemma 3.4, note that in this case the Betti numbers $\bar{\beta}_q=0$ for odd $q$, which together with Theorem 3.3  gives Morse type numbers $M_q=\bar{\beta}_q$, then $M_0=\bar{\beta}_0=1$ by (3.7) and (3.10), it contradicts the fact that $M_0\geq2$ and proves Theorem 1.2. \hfill\hb


\begin{thebibliography}{99}
\setlength{\itemsep}{-0.5mm}

\bibitem[Apo]{Apo} Tom M. Apostol, Modular functions and Dirichlet series in Number theory, GTM 41, 1990.
\bibitem[Ban]{Ban} V. Bangert, On the existence of closed geodesics on two-spheres.
{\it Inter. J. Math.} 4 (1993), 1-10.
\bibitem[BaL]{BaL} V. Bangert, Y. Long,   The existence of two closed
geodesics on every Finsler 2-sphere.  {\it Math. Ann.} 346 (2010), no. 2, 335-366.
%\bibitem[BCS]{BCS} D. Bao, S. S. Chern and Z. Shen, An Introduction to Riemann-Finsler Geometry. Springer. Berlin. 2000.
\bibitem[Bor]{Bor} A. Borel, Seminar on Transformation Groups. Princeton
Univ. Press. Princeton. 1960.
\bibitem[Bot]{Bot} Bott, R.,  On the iteration of closed geodesics and the Sturm
intersection theory.  {\it Comm. Pure Appl. Math.} 9 (1956) 171-206.
%\bibitem[BTZ]{BTZ} W. Ballmann, G. Thorbergsson and W. Ziller, Closed geodesics on positively curved manifolds. {\it Ann. of Math.} 116(1982),
213-247
% \bibitem[Cha]{Cha} K. C. Chang, Infinite Dimensional Morse Theory and Multiple Solution Problems. Birkh\"auser. Boston. 1993.
\bibitem[CoZ]{CoZ} C. Conley and E. Zehnder, Morse-type index theory for
flows and periodic solutions for Hamiltonian equations.
{\it Comm. Pure. Appl. Math.} 37 (1984) 207-253.
\bibitem[Dua]{Dua}H. Duan, Two elliptic closed geodesics on positively curved Finsler spheres. {\it J. Diff. Equ.}, 260 (2016), no. 12, 8388-8402.
\bibitem[DuL]{DuL1} H. Duan and Y. Long, Multiple closed geodesics on
bumpy Finsler $n$-spheres. {\it J. Diff. Equ.} 233 (2007), no. 1, 221-240.
\bibitem[DLW1]{DLW1} H. Duan, Y. Long and W. Wang, Two closed geodesics on compact
simply-connected bumpy Finsler manifolds. {\it J. Diff. Geom.} 104 (2016), no. 2, 275-289.
\bibitem[DLW2]{DLW2} H. Duan, Y. Long and W. Wang, The enhanced common index jump theorem for symplectic paths and non-hyperbolic closed geodesics on Finsler manifolds. {\it Calc. Var. PDE} 55 (2016), no. 6, Art. 145, 28 pp.
%\bibitem[FaR]{FaR} E. Fadell and P. Rabinowitz, Generalized cohomological index theories for Lie group actions with an application to bifurcation questions for Hamiltonian systems. {\it Invent. Math.} 45 (1978), no. 2, 139--174.
\bibitem[FHT]{FHT} Y.~F\'elix, S.~Halperin and J.-C.~Thomas, Rational Homotopy Theory,
Graduate Texts in Mathematics, vol.~205, Springer, 2001.
\bibitem[Fet]{Fet} A. Fet, A periodic problem in the calculus of variations. {\it Dokl. Akad. Nauk.
SSSR (N.S.) } 160 (1965) 287-289 [Russian]; {\it Soviet Math.}  6 (1965), 85-88.
\bibitem[Fra]{Fra} J. Franks,  Geodesics on $S^2$ and periodic points of
annulus diffeomorphisms. {\it Invent. Math.}  108 (1992), 403-418.
\bibitem[GG]{GG}V. Ginzburg, B. Gurel, Lusternik-Schnirelmann Theory and Closed Reeb Orbits,
{\it  Math. Z. } 295 (2020), 515-582.
\bibitem[GGM]{GGM}V. Ginzburg, B. Gurel and L. Macarini, Multiplicity of closed Reeb orbits on prequantization bundles. {\it Israel J. Math.} 228 (2018), no. 1, 407-453.
\bibitem[GrM]{GrM} D. Gromoll and W. Meyer, Periodic geodesics on compact Riemannian
manifolds. {\it J. Diff. Geom.}  3 (1969), 493-510.
\bibitem[GKV]{GKV} K.~Grove, S.~Halperin and M.~Vigu\'e-Poirrier, The rational homotopy theory of certain path spaces with
applications to geodesics, {\it Acta Math.} 140 (1978), 277-303.
\bibitem[GrR]{GrR} A. Granville and Z. Rudnick, Uniform distribution.
 In {\it  Equidistribution in Number Theory, An Introduction.}
(A. Granville and Z. Rudnick ed.) 1-13, (2007) Nato Sci. Series. Springer.
% \bibitem[HaW]{HaW} G. H. Hardy and E. M. Wright, An Introduction to the Theory of Numbers. 6th ed. Oxford University Press. 2008.
\bibitem[Hin]{Hin} N. Hingston,  On the growth of the number of closed
geodesics on the two-sphere. {\it Inter. Math. Research Notices.}
 9 (1993) 253-262.
 \bibitem[Kat]{Kat} A. B. Katok, Ergodic properties of degenerate
integrable Hamiltonian systems. {\it Izv. Akad. Nauk SSSR.} 37 (1973) (Russian),
{\it Math. USSR-Izv.} 7 (1973) 535-571.
\bibitem[Kli1]{Kli1} W. Klingenberg, Lectures on Closed Geodesics. Springer.
Berlin. 1978.
\bibitem[Kli2]{Kli2} W. Klingenberg,  Riemannian Geometry, Walter de Gruyter.
Berlin. 1st ed. 1982, 2nd ed, 1995.
\bibitem[Liu1]{Liu1} C. Liu,  The relation of the Morse index of closed geodesics
with the Maslov-type index of symplectic paths. {\it Acta Math. Sinica.} English
Series 21 (2005) 237-248.
\bibitem[Liu]{Liu} H. Liu, The optimal lower bound estimation of the number of closed
geodesics on Finsler compact space form $S^{2n+1}/\Gamma$,
{\it Calc. Var. PDE.} 58 (2019), no. 3, Art. 107, 21 pp.
% \bibitem[LLo]{LLo} C. Liu and Y. Long,  Iterated index formulae for closed geodesics with applications. {\it Science in China.} 45 (2002) 9-28.
\bibitem[LiX]{LiX} H. Liu and Y. Xiao, Resonance identity and
multiplicity of non-contractible closed geodesics on Finsler $\RP^n$.
{\it Adv. Math.} 318 (2017), 158-190.
\bibitem[LLX]{LLX} H. Liu, Y. Long and Y. Xiao, The existence of two non-contractible
closed geodesics on every bumpy Finsler compact space form,
{\it Discrete Contin. Dyn. Syst.} 38 (2018), no. 8, 3803-3829.
\bibitem[Lon1]{Lon1} Y. Long,  Bott formula of the Maslov-type index
theory. {\it Pacific J. Math.} 187 (1999), 113-149.
\bibitem[Lon2]{Lon2} Y. Long,  Precise iteration formulae of the
Maslov-type index theory and ellipticity of closed characteristics.
{\it Adv. Math.} 154 (2000), 76-131.
\bibitem[Lon3]{Lon3} Y. Long,  Index Theory for Symplectic Paths with
Applications. Progress in Math. 207, Birkh\"auser. Basel. 2002.
\bibitem[LoW]{LoW} Y. Long and W. Wang,  Stability of closed geodesics on Finsler
2-spheres. {\it J. Funct. Anal.} 255 (2008) 620-641.
\bibitem[LoZ]{LoZ} Y. Long and C. Zhu,  Closed characteristics on
compact convex hypersurfaces in $\R^{2n}$.  {\it Ann. of Math.}
155 (2002) 317-368.
\bibitem[LZe]{LZe} Y. Long and E. Zehnder, Morse theory for forced
oscillations of asymptotically linear Hamiltonian systems.
In {\it Stoc. Proc. Phys. and Geom.}, S. Albeverio et al. ed.
World Sci. (1990) 528-563.
\bibitem[LyF]{LyF} L. A. Lyusternik and A. I. Fet, Variational problems on
closed manifolds. {\it Dokl. Akad. Nauk SSSR (N.S.)} 81 (1951) 17-18
(in Russian).
% \bibitem[MaW]{MaW} J. Mawhin and M. Willem,  Critical Point Theory and Hamiltonian Systems. Springer. New York. 1989.
\bibitem[Rad1]{Rad1} H.-B. Rademacher,  On the average indices of closed
geodesics. {\it J. Diff. Geom.} 29 (1989), 65-83.
\bibitem[Rad2]{Rad2} H.-B. Rademacher, Morse Theorie und geschlossene
Geodatische. {\it Bonner Math. Schriften} Nr. 229 (1992).
\bibitem[Rad3]{Rad3} H.-B. Rademacher, The Fadell-Rabinowitc index and
closed geodesics. {\it J. London. Math. Soc.} 50 (1994) 609-624.
\bibitem[Rad4]{Rad4} H.-B. Rademacher, A Sphere Theorem for non-reversible
Finsler metrics. {\it Math. Ann.} 328 (2004) 373-387.
\bibitem[Rad5]{Rad5} H.-B. Rademacher, Existence of closed geodesics on
positively curved Finsler manifolds. {\it Ergodic Theory Dynam.
Systems.} 27 (2007), no. 3, 957--969.
 \bibitem[Rad6]{Rad6} H.-B. Rademacher, The second closed geodesic on Finsler
spheres of dimension $n>2$. {\it Trans. Amer. Math. Soc. } 362 (2010), no. 3, 1413-1421.
\bibitem[Rad7]{Rad7} H.-B. Rademacher, The second closed geodesic on the complex projective plane. {\it Front. Math. China.}  3 (2008), 253-258.
\bibitem[Sul]{Sullivan} D.~Sullivan, Infinitesimal computations in topology,
{\it Publ. Math. IH\'ES } 47(1977), 269-331.
\bibitem[Tai]{Taimanov} I.~A.~Taimanov, The spaces of non-contractible closed curves in compact
space forms, {\it Sb. Math.}  207:10 (2016), 1458-1470.
\bibitem[VPB]{VPB} M.~Vigu\'e-Poirrier and D.~Burghelea, A model for cyclic homology and algebraic $K$-theory of
1-connected topological spaces, {\it J. Diff. Geom.} 22 (1985), 243-253.
\bibitem[VPS]{VPS} M. Vigu\'e-Poirrier and D. Sullivan,  The homology theory of the closed
geodesic problem. {\it J. Diff. Geom.} 11 (1976), 633-644.
\bibitem[Wan1]{Wan1} W. Wang, Closed geodesics on positively curved Finsler spheres. {\it Adv.
Math.} 218 (2008) 1566-1603.
\bibitem[Wan2]{Wan2} W. Wang, On a conjecture of Anosov. {\it Adv. Math.} 230 (2012), 1597-1617.
\bibitem[Wan3]{Wan3} W. Wang,  On the average indices of closed geodesics on positively curved Finsler spheres. {\it Math. Ann.} 355 (2013), no. 3, 1049-1065.
%\bibitem[Wan4]{Wan4}W. Wang, Multiple closed geodesics on positively curved Finsler manifolds. {\it Adv. Nonlinear Stud. }19 (2019), no. 3, 495-518.
%\bibitem[Wil]{Wil} B. Wilking, Index parity of closed geodesics and rigidity of Hopf fibritions. {\it Invent. Math.} 144 (2001) 281-295.
\bibitem[Zil]{Zil} W. Ziller, Geometry of the Katok examples. {\it Ergod.
Th. \& Dynam. Sys.} 3 (1982) 135-157.

\end{thebibliography}
\end{document}